\documentclass[final]{siamart190516}

\usepackage{booktabs} 
\usepackage[noend]{algorithmic}
\newcommand{\STATEWITHCOMMENT}[2]{\STATE{\makebox[\widthof{bla bla bla bla bla bla bla}][l]{#1}$\triangleright$ #2}}
\usepackage{graphicx}
\usepackage{subcaption}
\usepackage{latexsym}
\usepackage{amsmath,amssymb} 
\usepackage{mathtools}       
\usepackage{stmaryrd}        
\usepackage{bm}              
\usepackage{hyperref}
\usepackage[utf8]{inputenc}
\usepackage{cite}
\usepackage{paralist}

\newsiamthm{assumption}{Assumption}
\crefname{assumption}{Assumption}{Assumptions}

\newcommand{\N}{\mathcal{N}}

\newcommand{\id}{\mathrm{I}}
\newcommand{\ind}[1]{1\left\{\scriptstyle{#1}\right\}}
\newcommand{\indlr}[1]{\ind{#1}}
\newcommand{\F}{\mathcal{F}}
\newcommand{\barsigma}{\bar{\sigma}}
\newcommand{\uu}{u}
\newcommand{\elll}{\ell}
\newcommand{\target}{\beta}
\newcommand{\ssf}{f_\mu}
\newcommand{\aup}{\alpha_{\uparrow}}
\newcommand{\adown}{\alpha_{\downarrow}}
\newcommand{\Ball}{\mathcal{B}}
\newcommand{\ptarg}{p_\mathrm{target}}
\newcommand{\Xab}{{\mathcal{X}_a^b}}

\providecommand{\argmin}{\operatornamewithlimits{argmin}} 
\providecommand{\limsup}{\operatornamewithlimits{limsup}} 
\providecommand{\liminf}{\operatornamewithlimits{liminf}} 
\DeclareMathOperator{\Cond}{Cond} 
\providecommand{\R}{\mathbb{R}} 
\providecommand{\E}{\mathbb{E}} 
\providecommand{\T}{\mathrm{T}} 
\renewcommand{\geq}{\geqslant} 
\renewcommand{\leq}{\leqslant} 
\DeclarePairedDelimiterX{\inner}[2]{\langle}{\rangle}{#1, #2}
\DeclarePairedDelimiter{\norm}{\lVert}{\rVert}
\DeclarePairedDelimiter{\abs}{\lvert}{\rvert}

\newcommand{\new}[1]{{\color{blue}#1}}
\newcommand{\del}[1]{{\color{orange}#1}}

\definecolor{tgcolor}{rgb}{0.8,0.2,0.2}

\newcommand{\correct}[2]{{\color{orange}#1}{\color{blue}#2}}

\renewcommand{\new}[1]{{#1}}
\renewcommand{\del}[1]{{}}
\renewcommand{\correct}[2]{#2}

\usepackage{enumitem}
\newlist{assump}{enumerate}{2}
\setlist[assump,1]{label=A\arabic*}

 \newsiamremark{remark}{Remark}

\title{Global Linear Convergence of Evolution Strategies on More than Smooth Strongly Convex Functions}

\author{Youhei Akimoto
	\thanks{Faculty of Engineering, Information and Systems, University of Tsukuba; RIKEN AIP, Tsukuba, Japan
	    (\email{akimoto@cs.tsukuba.ac.jp}).}
\and Anne Auger
	\thanks{Inria and CMAP, Ecole Polytechnique, IP Paris, France
		(\email{anne.auger@inria.fr}).}
\and Tobias Glasmachers
	\thanks{Institute for Neural Computation, Ruhr-University Bochum, Bochum, Germany
		(\email{tobias.glasmachers@ini.rub.de}).}
\and Daiki Morinaga
	\thanks{Department of Computer Science, University of Tsukuba; RIKEN AIP, Tsukuba, Japan
	    (\email{morinaga@bbo.cs.tsukuba.ac.jp}).}}

\headers{Global Linear Convergence of Evolution Strategies}{Y.~Akimoto, A.~Auger, T.~Glasmachers and D.~Morinaga}

\begin{document}

\maketitle

\begin{abstract}
Evolution strategies (ESs) are zeroth-order stochastic black-box optimization
heuristics invariant to monotonic transformations of the objective function. They evolve a multivariate normal distribution, from which candidate
solutions are generated. Among different variants, CMA-ES is nowadays
recognized as one of the state-of-the-art zeroth-order optimizers for
difficult problems. Albeit ample empirical evidence that ESs with a
step-size control mechanism converge linearly, theoretical guarantees
of linear convergence of ESs have been established only on limited
classes of functions. In particular, theoretical results on convex
functions are missing, where zeroth-order and also first-order
optimization methods are often analyzed.
In this paper, we establish almost sure linear convergence and a bound
on the expected hitting time of an \new{ES family, namely the $(1+1)_\kappa$-ES, which includes the (1+1)-ES with (generalized) one-fifth success rule} and an abstract covariance matrix adaptation with bounded condition number, on a broad class of functions.
The analysis holds for monotonic transformations of positively homogeneous functions and of quadratically
bounded functions,
the latter of which particularly includes monotonic transformation of strongly convex functions with
Lipschitz continuous gradient. As far as the authors know, this is the
first work that proves linear convergence of ES on such a broad class of
functions.
\end{abstract}





\providecommand{\ES}{(1+1)$\text{-ES}$}
\providecommand{\CMAES}{(1+1)$\text{-ES}_\kappa$}




\section{Introduction}

We consider the unconstrained minimization of an objective function $f: \R^d \to \R$ without the use of derivatives where an optimization solver sees $f$ as a zeroth-order {\it black-box oracle} \cite{nemirovski1995information, bubeck2014convex, nesterov2018lectures}. This setting is also referred to as derivative-free optimization \cite{ConnScheVice09}.
Such problems can be advantageously approached by randomized
algorithms that can typically be more robust to noise, non-convexity and
irregularities of the objective function than deterministic algorithms. There has been recently  a vivid interest in randomized derivative-free
algorithms giving rise to several theoretical studies of randomized direct search methods \cite{gratton2015direct},  trust region  \cite{bandeira2014convergence,gratton2017complexity} and model-based methods \cite{cartis2018global,paquette2018stochastic}.
We
refer to \cite{larson2019derivative} for an in-depth survey including
the references of this paragraph and additional ones.

In this context, we investigate Evolution Strategies (ES), which are among the oldest randomized derivative-free or zeroth-order black-box
methods \cite{Devroye:72,Schumer:Steiglitz:68,rechenberg:1973}. They are widely used in applications in different domains \cite{uhlendorf2012long,kriest2017calibrating,AAAI12-MacAlpine,geijtenbeek2013flexible,8080408,NIPS2018_7512,Volz:2018:EML:3205455.3205517,8698539,Dong_2019_CVPR,FUJII2018624}. Notably a specific ES called covariance-matrix-adaptation ES (CMA-ES) \cite{hansen:2001} is among the best solvers to address {\it difficult}
black-box problems. It is affine-invariant and implements complex adaptation
mechanisms for the sampling covariance matrix and step-size. It performs well on many ill-conditioned, non-convex,
non-smooth, and non-separable problems
\cite{Hansen:2010:CRA:1830761.1830790,Rios2013}.
ES are known to be difficult to analyze. Yet, given their importance in practice, it is essential to study them from a theoretical convergence perspective.

We focus on the arguably simplest and oldest adaptive
ES, denoted (1+1)-ES.
%
%
It samples a candidate solution from a Gaussian distribution whose
step-size (standard deviation) is adapted. The candidate solution is
accepted if and only if it is better than the current one (see
pseudo-code \Cref{algo}). The algorithm shares some similarities with
simplified direct search whose complexity analysis has been presented in
\cite{konevcny2014simple}. Yet the (1+1)-ES is comparison-based and thus
invariant to strictly increasing transformations of the objective
function. \new{Simplified direct search can be thought of as a variant of
mesh adaptive direct search \cite{audet2006mesh,abramson2009orthomads}.}
Arguably, in contrast to direct search, a sufficient decrease condition
cannot be guaranteed. This causes some difficulties for the analysis.
The (1+1)-ES is rotational invariant, while direct search candidate
solutions are created along a predefined set of vectors.
While the CMA-ES should always be preferred for practical applications over the (1+1)-ES variant analyzed here,
this latter variant  achieves faster linear convergence on well-conditioned
problems when compared to algorithms with established complexity
analysis (see \cite[Table~6.3 and Figure~6.1]{stich2013optimization} and
\cite[Figure~B.4]{auger2016linear} where the random pursuit algorithm
and the (1+1)-ES algorithms are compared, and also Appendix~\ref{app:NumericalResults}).

Prior theoretical studies of the \ES{} with $1/5$ success rule have established the global linear convergence on differentiable positively homogeneous functions (composed with a strictly increasing function) with a single optimum \cite{auger2016linear,AugerH13a}. Those results establish the almost sure linear convergence from all initial states. They however do not provide the dependency of the convergence rate with respect to the dimension. A more specific study on the sphere function $f(x) = \frac{1}{2} \|x\|^2$ establishes lower and upper bounds on the expected hitting time of an $\epsilon$-ball of the optimum in $\Theta(\log ( d \| m_0 - x^* \| /\epsilon))$, where $x^*$ is the optimum of the function, $m_0$ is the initial solution, and $d$ is the problem dimension \cite{akimoto2018drift}. Prior to that, a variant of the $(1+1)$-ES with one-fifth success rule had been analyzed on the sphere and certain convex quadratic functions establishing bounds on the expected hitting time \del{or }with overwhelming probability in $\Theta(\log ( \kappa_f d \| m_0 - x^* \| /\epsilon))$, where $\kappa_f$ is the condition number (the ratio between the greatest and smallest eigenvalues) of the Hessian \cite{jaegerskuepper2003analysis,jagerskupper2006quadratic,jagerskupper2007algorithmic,jagerskupper2005rigorous}.
%
Recently, the class of functions where the convergence of the \ES{} was
proven has been extended to continuously differentiable functions. This
analysis does not address the question of linear convergence, focusing
only on convergence as such, which is possibly
sublinear~\cite{Glasmachers2019ECJ}.

Our main contribution is as follows.
For a generalized version of the (1+1)-ES with one-fifth success rule, we prove bounds on the expected hitting time akin to linear convergence, i.e., hitting an $\epsilon$-ball in $\Theta(\log \| m_0 - x^*\|/\epsilon)$ iterations on a quite general class of functions. This class of functions includes all composites of Lipschitz-smooth strongly convex functions with a strictly increasing transformation. This latter transformation allows to include some non-continuous functions, and even functions with non-smooth level sets. We additionally deduce linear convergence with probability one.
Our analysis relies on finding an appropriate Lyapunov function with lower and upper-bounded expected drift. It is building on classical fundamental ideas presented by Hajek \cite{hajek1982hitting} and widely used to analyze stochastic hill-climbing  algorithms on discrete search spaces \cite{lengler2020drift}.

\subsection*{Notation}
Throughout the paper, we use the following notations.
The set of natural numbers $\{1,2,\ldots,\}$ is denoted $\mathbb{N}$.
Open, closed, and left open intervals on $\R$ are denoted by $( ,  )$, $[ ,  ]$, and $( ,  ]$, respectively.
The set of strictly positive real numbers is denoted by $\R_{>}$.
The Euclidean norm on $\R^d$ is denoted by $\norm{~}$.
Open and closed balls with center $c$ and radius $r$ are denoted as $\Ball(c, r) = \{x \in \R^d : \norm{x - c} < r\}$ and $\bar{\Ball}(c, r) = \{x \in \R^d : \norm{x - c} \leq r\}$, respectively.
Lebesgue measures on $\R$ and $\R^d$ are both denoted by the same symbol $\mu$. A multivariate normal distribution with mean $m$ and covariance matrix $\Sigma$ is denoted by $\N(m, \Sigma)$. Its probability measure and its induced probability density under Lebesgue measure are denoted by $\Phi(\cdot ; m, \Sigma)$
and $\varphi( \cdot; m, \Sigma)$. The indicator function of a set or condition $C$ is denoted by~$\ind{C}$.
\correct{We use Bachmann-Landau notations: $o(\cdot)$, $O(\cdot)$, $\Theta(\cdot)$, $\Omega(\cdot)$, $\omega(\cdot)$.}{%
We use Bachmann-Landau notations $f \in o(g)$, $O(g)$, $\omega(g)$, $\Omega(g)$ and $\Theta(g)$ to mean  $\limsup_{\epsilon \to 0} f(\epsilon) / g(\epsilon) = 0$,
$\limsup_{\epsilon \to 0} f(\epsilon) / g(\epsilon) \leq C$ for some $C > 0$,
$\liminf_{\epsilon \to 0} f(\epsilon) / g(\epsilon) = \infty$,
$\liminf_{\epsilon \to 0} f(\epsilon) / g(\epsilon) \geq C$ for some $C > 0$,
and $f \in O(g)$ and $f \in \Omega(g)$, respectively.}

\section{Algorithm, Definitions and Objective Function Assumptions}
\label{sec:algorithm_and_definition}

%
%

\subsection{Algorithm: \ES{}  with Success-based Step-size Control}\label{sub:es}

We analyze a generalized version of the (1+1)-ES with one-fifth success rule 
presented in \Cref{algo}, which implements one of the oldest approaches to adapt the step-size in randomized optimization methods \cite{rechenberg:1973,Devroye:72,Schumer:Steiglitz:68}. The specific implementation was proposed in \cite{kern:2004}.
At each iteration, a candidate solution $x_t$ is sampled. It is centered in the current incumbent $m_t$ and follows a multivariate normal distribution with mean vector $m_t$ and covariance matrix equal to $\sigma_t^2 I_d$ where $I_d$ denotes the identity matrix. The candidate solution is accepted, that is $m_t$ becomes $x_t$, if and only if $x_t$ is better than $m_t$ (i.e.\ $f(x_t) \leq f(m_t)$). In this case, we say that the candidate solution is successful. The step-size $\sigma_t$ is adapted so as to maintain a probability of success
to be approximately the target success probability denoted by
$\ptarg := \frac{\log(1 / \adown)}{ \log(\aup / \adown)}$.
To do so, the step-size is increased by the increase factor $\aup > 1$ in case of success (which is an indication that the step-size is likely to be too small) and decreased by the decrease factor $\adown < 1$ otherwise.
The covariance matrix $\Sigma_t$ of the sampling distribution of candidate solutions is adapted in the set $\mathcal{S}_\kappa$ of positive-definite symmetric matrices with determinant $\det(\Sigma) = 1$ and condition number $\Cond(\Sigma) \leq \kappa$.
We do not assume any specific update mechanism for $\Sigma$, but we assume that the update of $\Sigma$ is invariant to any strictly increasing transformation of $f$. We call such an update comparison-based (see Lines~7 and 11 of \Cref{algo}).
 Then, our algorithm behaves exact-equally on $f$ and on $g \circ f$ for all strictly increasing functions $g:\R \to \R$ \big(i.e., $g(s) \lesseqgtr g(t) \Leftrightarrow s \lesseqgtr t$\big). This defines a class of comparison-based randomized algorithms and we denote it as \CMAES.
For $\kappa = 1$, it is simply denoted as \ES.

\begin{algorithm}
\caption{\CMAES{} with success-based step-size adaptation}
\label{algo}
\begin{algorithmic}[1]
\STATE{\textbf{input} $m_0 \in \mathbb{R}^d$, $\sigma_0 > 0$, $\Sigma_0 = I$, $f: \R^d \to \R$}, \textbf{parameter} $\aup > 1 > \adown > 0$
\FOR {$t = 1,2,\dots$, \textit{until stopping criterion is met}}
	\STATE {sample $x_t \sim m_t + \sigma_t   \N(0,  \Sigma_t)$}
	\IF {$f\big(x_t\big) \leq f\big(m_t\big)$}
		\STATEWITHCOMMENT{$m_{t+1} \leftarrow x_t$}{move to the better solution}
		\STATEWITHCOMMENT{$\sigma_{t+1} \leftarrow \sigma_t   \aup$}{increase the step size}
\STATEWITHCOMMENT{$\Sigma_{t+1} \in \mathcal{S}_\kappa$}{adapt the covariance matrix}
	\ELSE
		\STATEWITHCOMMENT{$m_{t+1} \leftarrow m_t$}{stay where we are}
		\STATEWITHCOMMENT{$\sigma_{t+1} \leftarrow \sigma_t   \adown$}{decrease the step size}
\STATEWITHCOMMENT{$\Sigma_{t+1} \in \mathcal{S}_\kappa$}{adapt the covariance matrix}
	\ENDIF
\ENDFOR
\end{algorithmic}
\end{algorithm}
Note that $\aup$ and $\adown$ are not meant to be tuned depending on the function properties. How to choose such constants for $\Sigma_t = I_d$ is well-known and is related to the so-called evolution window \cite{rechenberg1994evolutionsstrategie}.
In practice, $\adown = \aup^{-1/4}$ is the most commonly used setting, which leads to $\ptarg = 1/5$. It has been shown to be close to optimal, which gives nearly optimal (linear) convergence rate on the sphere function \cite{rechenberg:1973,Devroye:72}.
Hereunder we write $\theta = (m, \sigma, \Sigma)$ as the state of the algorithm, $\theta_t =(m_t,\sigma_t,\Sigma_t)$ and the state-space is denoted by $\Theta$.

\begin{figure}
\includegraphics[trim=0.3cm 0cm 0.3cm 0.4cm,height=3.2cm]{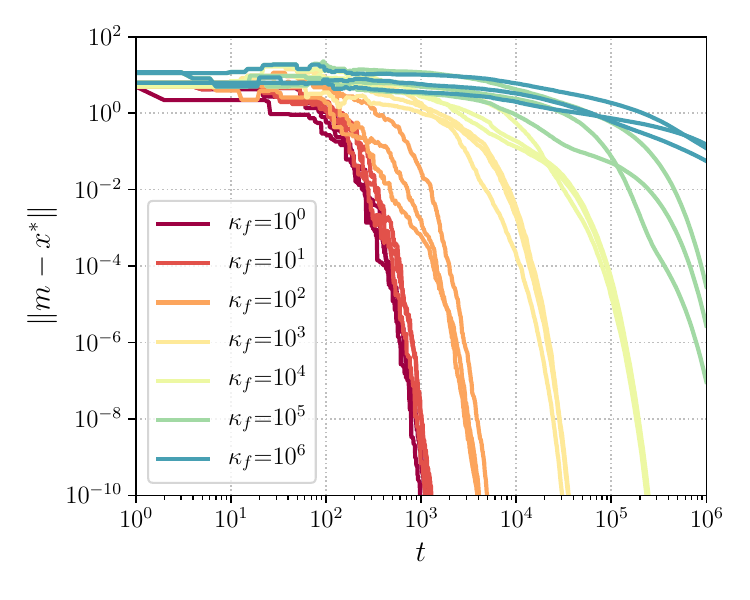}%
\includegraphics[trim=0.3cm 0cm 0.3cm 0.4cm,height=3.2cm]{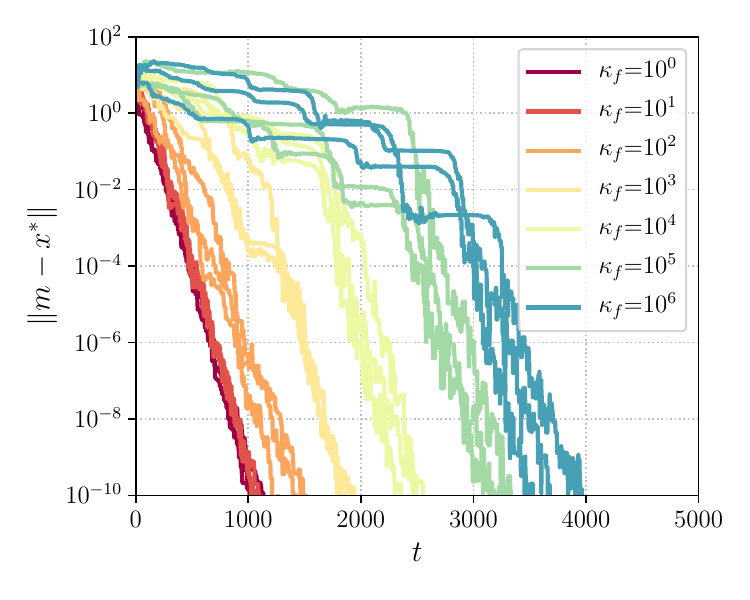}%
\includegraphics[trim=0.3cm -0.52cm 0.8cm 1cm,clip,height=3.2cm]{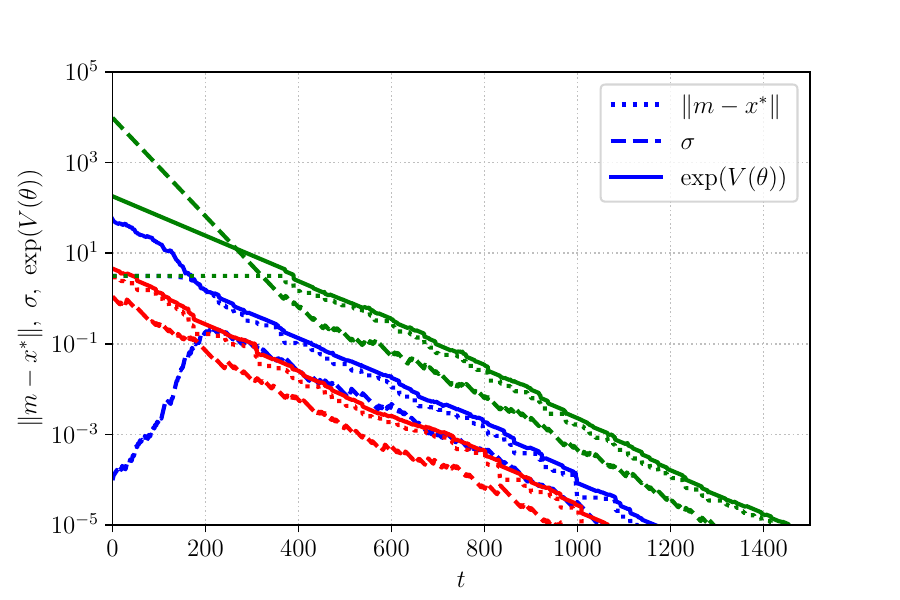}
%
\caption{\new{Convergence of the (1+1)-ES (left) and the (1+1)-CMA-ES (middle) on $10$ dimensional ellipsoidal function $f(x) = \frac12\sum_{i=1}^{d} \kappa_f^{\frac{i-1}{d-1}} x_i^2$ with $\kappa_f = 10^0, 10^1, \dots, 10^6$. The y-axis displays the distance to the optimum (and not the function value).
We employ the covariance matrix adaptation mechanism proposed by \cite{Arnold:2010:ACM:1830483.1830556}, where $\sigma$ is adapted as in \Cref{algo} with $\aup=e^{0.1}$ and $\adown=e^{-0.025}$. Note the logarithmic scale of the time axis of the left plot vs.\ the linear time axis of the middle plot.\\
Right: Three runs of \ES{} ($\aup = e^{0.1}$ and $\aup = e^{-0.025}$) on $10$ dimensional spherical function $f(x) = \frac12\norm{x - x^*}^2$ with initial step-size $\sigma_0 = 10^{-4}$, $1$, and $10^{4}$ (in blue, red, green, respectively). Plotted are the distance to the optimum (dotted line), the step-size (dashed line), and the potential function $V(\theta)$ defined in \eqref{eq:potentialfunction} (solid line) with $v = 4/d$, $\ell = \aup^{-10}$, and $u = \adown^{-10}$.
}
}
\label{fig:on_quad}
\end{figure}

\Cref{fig:on_quad} shows typical runs of the (1+1)-ES and a version of \CMAES{} proposed in \cite{Arnold:2010:ACM:1830483.1830556}, which is known as the (1+1)-CMA-ES, on a $10$-dimensional ellipsoidal function with different condition numbers $\kappa_f$ of the Hessian.
It is empirically observed that $\Sigma_t$ in the (1+1)-CMA-ES approaches the inverse Hessian $\nabla^2 f(m_t)$ of the objective function up to the scalar factor if the objective function is convex quadratic.
The runtime of (1+1)-ES scales linearly with $\kappa_f$ (notice the logarithmic scale of the horizontal axis), while the runtime of the (1+1)-CMA-ES suffers only an additive penalty,
roughly proportional to the logarithm of $\kappa_f$. Once the Hessian is well approximated by $\Sigma$ (up to a scalar factor), it approaches the global optimum geometrically at the same rate for different values of $\kappa_f$.

In our analysis, we do not assume any specific $\Sigma$ update mechanism, hence it does not necessarily behave as shown in \Cref{fig:on_quad}. Our analysis is therefore the worst case analysis (for the upper bound of the runtime) and the best case analysis (for the lower bound of the runtime) among the algorithms in \CMAES.

\subsection{Preliminary Definitions}

\subsubsection{Spatial Suboptimality Function}\label{sub:ssf}

The algorithms studied in this paper are comparison-based
 and thus invariant to strictly increasing transformations of $f$. If the convergence of the algorithms is measured in terms of $f$, say by investigating the convergence or hitting time of the sequence $f(m_t)$, this will not reflect the invariance to monotonic transformations of $f$ because the first iteration $t_0$ such that $f(m_{t_0})\leq \epsilon$ is not equal to the first iteration $t_0'$ such that $g(f(m_{t_0'})) \leq \epsilon$ for some $\epsilon > 0$. For this reason, we introduce a quality measure called \emph{spatial suboptimality function} \cite{Glasmachers2019ECJ}. It is the $d$th root of the volume of the sub-levelset where the function value is better or equal to $f(x)$:
\begin{definition}[Spatial Suboptimality Function]\label{def:ssf}
Let $f:\R^d \to \R$ be a measurable function with respect to the Borel $\sigma$ algebra of $\R^d$ (simply referred to as measurable function in the sequel). Then the spatial suboptimality function $\ssf:\R^d \to [0, +\infty]$ is defined as
\begin{equation}
\ssf(x) = \sqrt[d]{\mu\left(f^{-1}\left((-\infty,f(x)]\right)\right)}  = \sqrt[d]{\mu\left(\left\{y \in \R^d \,\big|\, f(y) \leq f(x)\right\}\right)} \label{eq:ssf} \enspace.
\end{equation}
\end{definition}
We remark that for any $f$, the suboptimality function $f_\mu$ is greater or equal to zero. For any $f$ and any strictly increasing function $g:{\rm Im}(f) \to \R$, $f$ and its composite $g \circ f$ have the same spatial suboptimality function such that hitting time of $f_\mu$ smaller than $\epsilon > 0$ will be the same for $f$ or $g \circ f$. 
Moreover, there exists a strictly increasing function $g$ such that $\ssf(x) = g(f(x))$ holds $\mu$-almost everywhere \cite[Lemma~1]{Glasmachers2019ECJ}.

We will investigate the expected first hitting time of $\norm{m_t - x^*}$ to $\epsilon > 0$. For this, we will bound the
 first hitting time of $\norm{m_t - x^*}$ to $\epsilon$ by the first hitting time of $f_\mu(m_t)$ to a constant times $\epsilon$.
 To understand why, consider  first a strictly convex quadratic function $f$ with Hessian $H$ and minimal solution $x^*$.
 We have $\ssf(x) = V_d   \big[ 2 (f(x)-f(x^*)) /  \det(H)^{1/d} \big]^{1/2}$ for all $x \in \R^d$, where $V_d = \pi^{1/2} / \Gamma^{1/d}(d/2 + 1)$ is the $d$th root of the volume of the $d$-dimensional unit hyper-sphere \cite{Akimoto2012gecco}. This implies that the first hitting time of $\ssf(m_t)$ translates to the first hitting time of
 $\sqrt{f(m_t) - f(x^*)}$.
 \new{We have} $\sqrt{\lambda_{\min}} \norm{x - x^*} \leq \sqrt{f(x) - f(x^*)} \leq \sqrt{\lambda_{\max}} \norm{x - x^*}$\new{, where $\lambda_{\min}$ and $\lambda_{\max}$ are the minimal and maximal eigenvalues of $H$}.
 E.g., consider $f(x) = \norm{x - x^*}^2 + 1$. \new{Then the above condition} also translates to the first hitting time of $\norm{m_t - x^*}$.
More generally, we will formalize an assumption on $f$ later on (Assumption~\ref{asm:radius}), which allows us to bound $\norm{x - x^*}$ by a constant times $\ssf(x)$ from above and below (see \eqref{eq:kind-of-quad-bound}), implying that the first hitting time of $\norm{m_t - x^*}$ to $\epsilon$ is bounded by that of $\ssf(m_t)$ to $\epsilon$, times a constant.

\subsubsection{Success Probability}\label{sub:sp}

The success probability, i.e., the probability of sampling a candidate solution $x_t$ with an objective function better than or equal to that of the current solution $m_t$, plays an important role in the analysis of the \CMAES{} with success-based step-size control mechanism. We present here several useful definitions related to the success probability.

We start with the definition of the \emph{success domain with rate $r$}
and the \emph{success probability with rate $r$}.
The probability to sample in the $r$-success domain is called success probability with rate $r$. When $r=0$ we simply talk about success probability.%
\footnote{For $r = 0$,
the success domain $S_0(m)$ is not necessarily equivalent to the
sub-levelset $S_0'(m) := \{x \in \R^d \mid f(x) \leq f(m) \}$, where it
always holds that $S_0'(m) \subseteq S_0(m)$. However, since it is
guaranteed that $\mu(S_0(m) \setminus S_0'(m)) = 0$ by
\cite[Lemma~1]{Glasmachers2019ECJ}, due to the absolute continuity of
$\Phi( ; 0, \Sigma)$ for $\Sigma \in \mathcal{S}_\kappa$, the
success probability with rate $r = 0$ is equivalent to $\Pr_{z \sim
\mathcal{N}(0, \Sigma)}\left[ m + \ssf(m) \cdot\bar\sigma z \in S_0'(m)
\right]$, with $\bar\sigma$ defined in \eqref{eq:barsigma}.}
\begin{definition}[Success Domain]\label{def:sd}
For a measurable function $f:\R^d \to \R$ and $m \in \R^d$ such that $\ssf(m) < \infty$, the $r$-success domain at $m$ with  $r \in [0, 1]$ is defined as
\begin{equation}
S_r(m) = \{ x \in \R^d \mid \ssf(x) \leq (1 - r) \ssf(m)\} \enspace.
\end{equation}
\end{definition}

\begin{definition}[Success Probability]\label{def:progress}
Let $f$ be a measurable function and let $m_0 \in \R^d$ be the initial search point satisfying $\ssf(m_0) < \infty$. For any $r \in [0, 1]$ and any $m \in S_0(m_0)$, the success probability with rate $r$ at $m$ under the normalized step-size $\bar\sigma$
	is defined as
	\begin{equation}
	\label{eq:barsigma}
	p^{\mathrm{succ}}_{r}(\bar\sigma; m,\Sigma) = \Pr_{z \sim \mathcal{N}(0, \Sigma)}\left[  m + \ssf(m)\bar\sigma z \in S_r(m) \right] \enspace.
	\end{equation}
\end{definition}



\Cref{def:progress} introduces the notion of \emph{normalized step-size $\bar\sigma$} and the success probability is defined as a function of $\barsigma$ rather than the actual step-size $\sigma=\ssf(m)   \bar\sigma$. This is motivated by the fact that as $m$ approaches the global optimum $x^*$ of $f$, the step-size $\sigma$ needs to shrink for the success probability to be constant. If the objective function is $f(x) = \frac12 \norm{x-x^*}^2$ and the covariance matrix is the identity matrix, then the success probability is fully controlled by
$\barsigma_t = \sigma_t / \ssf(m_t) \propto \sigma_t / \norm{m_t-x^*}$ and is independent of $m_t$. This statement can be formalized in the following way.
\begin{lemma}\label{lemma:r-proba-success-sphere}
If $f(x) = \frac12 \norm{x-x^*}^2$, then letting $e_1=(1,0,\ldots,0)$, we have
$$
p^{\mathrm{succ}}_{r}(\bar\sigma; m,\id) = \Pr_{z \sim \mathcal{N}(0, \id)} \left[ m +  \ssf(m) \bar\sigma z \in S_r(m) \right] =  \Pr_{z \sim \mathcal{N}(0, \id)} \left[  \| e_1 +  V_d \bar\sigma  z \| \leq (1-r) \right] \enspace.
$$%
\end{lemma}
\begin{proof}
The suboptimality function is the $d$-th rooth of the volume of a sphere of radius $\| x - x^*\|$. Hence $f_\mu(x) = V_d \norm{ x - x^* }$. Then, the proof follows the derivation in Section~3 in \cite{akimoto2018drift}.
\end{proof}%
Therefore, $\barsigma$ is more discriminative than $\sigma$ itself. In general, the optimal step-size is not necessarily proportional to neither $\norm{m_t-x^*}$ nor $\ssf(m_t)$.

Since the success probability under a given normalized step-size depends on $m$ and $\Sigma$, we define the upper and lower success probability as follows.
\begin{definition}[Lower and Upper Success Probability]\label{def:lusp}
Let $\mathcal{X}_{a}^{b} = \{x \in \R^d : a < \ssf(x) \leq b\}$.
	Given the normalized step-size $\bar\sigma > 0$, the lower and upper success probabilities are defined as
	\begin{align*}
	p^\mathrm{lower}_{(a, b]}(\bar\sigma) &= \inf_{m \in \mathcal{X}_{a}^{b}} \inf_{\Sigma \in \mathcal{S}_\kappa} p^{\mathrm{succ}}_{0}(\bar\sigma; m,\Sigma) \enspace,& 
	p^\mathrm{upper}_{(a, b]}(\bar\sigma) &= \sup_{m \in \mathcal{X}_{a}^{b}} \sup_{\Sigma \in \mathcal{S}_\kappa} p^{\mathrm{succ}}_{0}(\bar\sigma; m,\Sigma) \enspace.
	\end{align*}
\end{definition}

A central quantity for our analysis is the limit for $\bar\sigma$ to $0$ of the success probability $p^\mathrm{succ}_0 (\bar\sigma; m, \Sigma)$. Intuitively, if this limit is too small for a given $m$ (compared to $p_{\rm target}$), because the ruling principle of the algorithm is to decrease the step-size if the probability of success is smaller than $p_{\rm target}$,
the step-size will keep decreasing, causing undesired convergence.
Following Glasmachers~\cite{Glasmachers2019ECJ}, we introduce the concepts of \emph{$p$-improbability} and \emph{$p$-criticality}. They are defined in \cite{Glasmachers2019ECJ} by the probability of sampling a better point from the isotropic normal distribution in the limit of the step-size to zero. Here, we define $p$-improvability and $p$-criticality for a general multivariate normal distribution.

\begin{definition}[$p$-improvability and $p$-criticality]
Let $f$ be a measurable function. The function $f$ is called $p$-improvable at $m \in \R^d$ under the covariance matrix $\Sigma \in \mathcal{S}_\kappa$ if there exists $p \in (0, 1]$ such that
\begin{equation}
p = \liminf_{\bar\sigma \to +0} p^\mathrm{succ}_0 (\bar\sigma; m, \Sigma) \enspace.
\end{equation}
Otherwise, it is called $p$-critical.
\end{definition}

The connection to the classical definition of the critical points for continuously differentiable functions is summarized in the following proposition, which is an extension of Lemma~4 in \cite{Glasmachers2019ECJ}, taking a non-identity covariance matrix into account.
\begin{proposition}\label{prop:asm2}

Let $f = g \circ h$ be a measurable function where $g$ is any strictly
increasing function and $h$ is continuously differentiable. Then, $f$ is
$p$-improvable with $p = 1/2$ at any regular point $m$ where
$\nabla h(m) \neq 0$ under any $\Sigma \in \mathcal{S}_\kappa$.
Moreover, if $h$ is twice continuously differentiable at a critical
point $m$ where $\nabla h(m) = 0$ and at least one eigenvalue of
$\nabla^2 f(m)$ is non-zero, under any $\Sigma \in \mathcal{S}_\kappa$,
$m$ is $p$-improvable with $p = 1$ if $\nabla^2 h(m)$ has only
non-positive eigenvalues, $p$-critical if $\nabla^2 h(m)$ has only
non-negative eigenvalues, and $p$-improvable with
some $p > 0$ if $\nabla^2 h(x)$ has at least one strictly negative eigenvalue.
\end{proposition}
\begin{proof}
  Note that $p^\mathrm{succ}_0 (\bar\sigma; m, \Sigma)$ on $f$ is equivalent to $p^\mathrm{succ}_0 (\bar\sigma; m, I_d)$ on $f'(x) = f(m + \sqrt{\Sigma} (x - m))$. Therefore, it suffices to show that the claims hold for $\Sigma = I_d$ on $f'$, which is proven in Lemma~4 in \cite{Glasmachers2019ECJ}.
\end{proof}%

\subsection{Main Assumptions on the Objective Functions}
\label{sec:assump}

\begin{figure}[t]
	\centering
	\includegraphics[height=4.8cm]{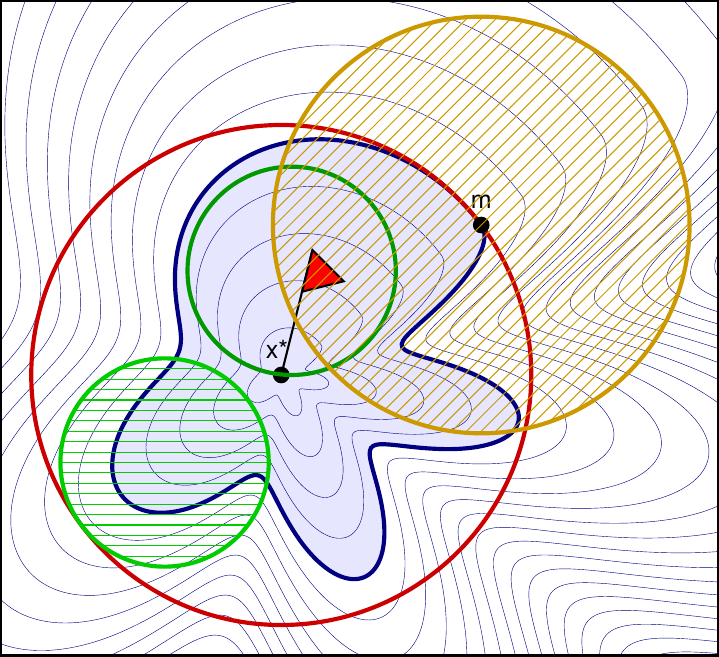}\
	\includegraphics[height=4.8cm]{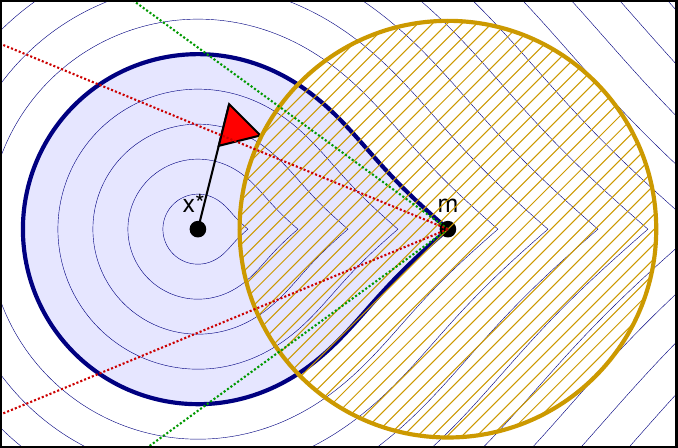}%
	\caption{
		\label{fig:levelset-balls}
		The sampling distribution is indicated by the mean $m$ and the
		shaded orange circle, indicating one standard deviation. The
		blue set is the sub-levelset $S_0(m)$ of points improving upon
		$m$.
		{\bf Left:}
        Illustration of property~\ref{asm:radius} in \Cref{sec:assump}.
		 The blue set is enclosed in the red (outer) ball of radius
		$C_u \ssf(m)$ and contains the dark green (inner) ball of radius
		$C_{\ell} \ssf(m)$. The shaded light green ball indicates the
		worst case situation captured by the bound, namely that the
		small ball is positioned within the large ball at maximal
		distance to~$m$.
		{\bf Right:}
        Illustration of property~\ref{asm:levelset} in \Cref{sec:assump}.
        Although the level set has a kink at $m$, there exists a
		cone centered at $m$ covering a probability mass of
		$p^\mathrm{limit}$ of improving steps (inside $S_0(m)$) for
		small enough step size $\sigma$ (green outline). It contains a
		smaller cone (red outline) covering a probability mass
		of~$p^\mathrm{target}$.
	}
\end{figure}

Given positive real numbers $a$ and $b$ satisfying $0 \leq a < b \leq + \infty$, and a measurable objective function,
\new{let $\Xab$ be the set defined in \Cref{def:lusp}}.



We pose two core assumptions on the objective functions under which we will derive an upper bound on the expected first hitting time of $[0,\epsilon]$ by $\ssf(m_t)$ (\Cref{theo:FHT-UB}) provided $a \leq \epsilon \leq \ssf(m_0) \leq b$. First, we require to be able to embed and include balls of radius scaling with $\ssf(m)$ into the sublevel sets of $f$. We do not require this to hold on the whole search space but, for a set $\mathcal{X}_a^b$. 

\begin{assump}[resume]
\item\label{asm:radius} We assume that $f$ is a measurable function and that there exists
  $\mathcal{X}_a^b$ such that for any $m \in \mathcal{X}_a^b$, there exist an open ball $\Ball_\ell$ with radius $C_{\ell} \ssf(m)$ and a closed ball $\bar{\Ball}_u$ with radius $C_u \ssf(m)$ such that it holds $\Ball_{\ell} \subseteq \{x \in \R^d \mid \ssf(x) < \ssf(m) \}$ and
  $\{x \in \R^d \mid \ssf(x) \leq \ssf(m) \} \subseteq \bar{\Ball}_u$.
\end{assump}
We do not specify the center of those balls that may or may not be centered on an optimum of the function.
We will see in \Cref{pro:pupper} that this assumption allows to bound $p^\mathrm{lower}_{(a, b]}(\bar\sigma)$ and $p^\mathrm{upper}_{(a, b]}(\bar\sigma)$ by tractable functions of $\bar\sigma$ which will be essential for the analysis.
The property is illustrated in \Cref{fig:levelset-balls}.

 The second assumption requires that the functions are $p$-improvable for $p$ which is lower-bounded uniformly over $\Xab$. 
\begin{assump}[resume]
\item\label{asm:levelset} Let $f$ be a measurable function, we assume that there exists $\Xab$ and there exists $p^\mathrm{limit} > p^\mathrm{target}$ such that for any $m \in \mathcal{X}_{a}^{b}$ and any $\Sigma \in \mathcal{S}_\kappa$, the objective function $f$ is $p$-improvable for some $p \geq p^\mathrm{limit}$, i.e.,
	\begin{equation}
		\liminf_{\bar\sigma \downarrow 0}  p^\mathrm{lower}_{(a, b]}(\bar\sigma) \geq p^\mathrm{limit} \enspace.
		\label{eq:asm:levelset}
	\end{equation}
      \end{assump}
	The property is illustrated in \Cref{fig:levelset-balls}.
This assumption implies in particular for a continuous function that $\Xab$ does not contain any local optimum.
This latter assumption is required to obtain global convergence \cite[Theorem~2]{Glasmachers2019ECJ} even without any covariance matrix adaptation (i.e.\ with $\kappa = 1$) and it can be intuitively understood:
If we have a point which is $p$-improvable with $p < \ptarg$ and which is not a local minimum of the function, then, starting with a small step-size, the success-based step-size control may keep decreasing the step-size at such a point and the \CMAES{} will prematurely converge to a point that is not a local optimum.


If \ref{asm:radius} is satisfied with balls centered at the optimum $x^*$ of the function $f$, then it is easy to see that for all $x \in \Xab$
\begin{equation}\label{eq:kind-of-quad-bound}
C_\ell \ssf(x) \leq \norm{x - x^*} \leq C_u \ssf(x) \enspace.
\end{equation}
If the balls are not centered at the optimum, we have the one-side inequality $\norm{x - x^*} \leq 2 C_u \ssf(x)$. Hence, the expected first hitting time of $\ssf(m_t)$ to $[0, \epsilon]$ translates to an upper bound for the expected first hitting time of $\norm{m_t - x^*}$ to $[0, 2 C_u \epsilon]$.


We remark that \ref{asm:radius} and \ref{asm:levelset} satisfied for $a=0$ allow to include \new{some non-differentiable} functions with non-convex sublevel sets as illustrated in \Cref{fig:levelset-balls}.

  We now give two examples of functions that satisfy \ref{asm:radius} and \ref{asm:levelset}, including function classes where linear convergence of numerical optimization algorithms are typically analyzed.
The first class consists of quadratically bounded functions. It includes all strongly-convex functions with Lipschitz continuous gradient. It also includes some non-convex functions.
The second class consists of positively homogeneous functions.
The levelsets of a positively homogeneous function are all geometrically similar around $x^*$.
\begin{assump}[resume]
\item\label{quad-bound} We assume that $f = g \circ h$ where $g$ is a strictly increasing function and $h$ is measurable, continuously differentiable with the unique critical point $x^*$, and quadratically bounded around $x^*$, i.e., for some $L_{u} \geq L_{\ell} > 0$,
\begin{align}
  (L_{\ell} / 2) \norm{x-x^*}^2 \leq h(x) - h(x^*) \leq ( L_u / 2) \norm{x-x^*}^2
  \enspace.
\label{eq:quad-bound}
\end{align}
\item\label{asm:pos-hom} We assume that $f = g \circ h$ where $h$ is continuously differentiable and positively homogeneous with a unique optimum $x^*$, i.e., for some $\gamma > 0$
\begin{equation}\label{eq:pos-homogeneous}
  h(x^* + \gamma  x) = h(x^*) + \gamma  \left( h(x^* + x) - h(x^*) \right)
  \enspace.
\end{equation}
\end{assump}

The following lemmas show that these assumptions imply \ref{asm:radius} and \ref{asm:levelset}.
The proofs of the lemmas are presented in \Cref{apdx:lem:asm1} and \Cref{apdx:prop:pos-hom}, respectively.
\begin{lemma}\label{lem:asm1}
Let $f$ satisfy \ref{quad-bound}. Then, $f$ satisfies \ref{asm:radius} and \ref{asm:levelset} with $a=0$, $b = \infty$,
$C_{\ell} = \frac{1}{V_d}\sqrt{\frac{L_{\ell}}{L_u}}$ and $C_u = \frac{1}{V_d}\sqrt{\frac{L_u}{L_{\ell}}}$.
\end{lemma}
\begin{lemma}\label{prop:pos-hom}
Let $f$ be positively homogeneous satisfying \ref{asm:pos-hom}, then the suboptimality function $\ssf(x)$ is proportional to $h(x) - h(x^*)$ and satisfies \ref{asm:radius} and \ref{asm:levelset} for $a=0$ and $b=\infty$ with $C_u = \sup \{\norm{x - x^*} : \ssf(x) = 1 \}$ and $C_{\ell} = \inf \{\norm{x - x^*} : \ssf(x) = 1 \}$.
\end{lemma}%

\section{Methodology: Additive Drift on Unbounded Continuous Domains}
\label{section:drift}


\subsection{First Hitting Time}\label{sub:fht}
We start with the generic definition of the \emph{first hitting time} of a stochastic process $\{X_t : t \geq 0\}$, defined as follows.
\begin{definition}[First hitting time]\label{def:firsthittingtime}
Let $\{X_t : t \geq 0\}$ be a sequence of real-valued random variables adapted to the natural filtration $\{\mathcal{F}_t : t \geq 0\}$ with initial condition $X_0 = \target_0 \in \R$. For $\target < \target_0$, the first hitting time $T_{\target}^X$ of $X_t$ to the set $(-\infty, \target]$ is defined as $T_{\target}^X = \inf\{t : X_t \leq \target\}$.
\end{definition}

The first hitting time is the number of iterations that the stochastic process requires to reach the target level $\target < \target_0$ for the first time. In our situation, $X_t = \norm{m_t - x^*}$ measures the distance from the current solution $m_t$ to the target point $x^*$ (typically, global or local optimal point) after $t$ iterations. Then, $\target = \epsilon > 0$ defines the target accuracy and $T_{\epsilon}^X$ is the runtime of the algorithm until it finds an $\epsilon$-neighborhood $\Ball(x^*, \epsilon)$. The first hitting time $T_{\epsilon}^X$ is a random variable as $m_t$ is a random variable. In this paper, we focus on the \emph{expected first hitting time} $\E[T_{\epsilon}^X]$. We want to derive lower and upper bounds on this expected hitting time that relate to the linear convergence of $X_t$ towards $x^*$. Such bounds take the following form:
There exist $C_T, \tilde{C}_T \in \R$ and $C_R > 0$, $\tilde{C}_R >0 $ such that for any $0 < \epsilon \leq \beta_0$
\begin{equation}\label{eq:FHT_upper}
\tilde{C}_T +  \frac{\log\left({\norm{m_0 - x^*}}/{\epsilon}\right)}{\tilde{C}_R} \leq \E[T_{\epsilon}^X | \mathcal{F}_0] \leq C_T +  \frac{\log(\norm{m_0 - x^*} / \epsilon)}{C_R} \enspace.
\end{equation}
That is, the time to reach the target accuracy scales logarithmically with the ratio between the initial accuracy $\norm{m_0 - x^*}$ and the target accuracy $\epsilon$. The first pair of constants, $C_T$ and $\tilde{C}_T$, capture the transient time, which is the time that adaptive algorithms typically spend for adaptation. The second pair of constants, $C_R$ and $\tilde{C}_R$, reflect the speed of convergence (logarithmic convergence rate). Intuitively, assuming that $C_R$ and $\tilde{C}_R$ are close, the distance to the optimum decreases in each step at a rate of approximately $\exp(-C_R) \approx \exp(-\tilde{C}_R)$.
While upper-bounds inform us about the (linear) convergence, the lower-bound helps understanding whether the upper bounds are tight.

Alternatively,
linear convergence can be defined as the following property: there exits $C > 0$ such that
\begin{equation}\label{eq:asymptot-lin-conv}
\limsup_{t \to \infty}\frac{1}{t} \log\frac{\|m_t - x^* \|}{\|m_0 - x^*\|} \leq - C \,\,  {\rm almost\,\, surely.}
\end{equation}
When we have an equality in the previous statement, we say that $\exp(-C)$ is the convergence rate.

\Cref{fig:on_quad} (right plot) visualizes three different runs of the \ES{} on a function with spherical level sets with different initial step-sizes. First of all, we clearly observe linear convergence. The first hitting time of $\Ball(x^*, \epsilon)$ scales linearly with $\log(1/\epsilon)$ for a sufficiently small $\epsilon > 0$. Second, its convergence speed is independent of the initial condition. Therefore, we expect to have  universal constants $C_R$ and $\tilde{C}_R$ independent of the initial state. Last, depending on the initial step-size, the transient time can vary. If the initial step-size is too large or too small, it does not produce progress in terms of $\norm{m_t - x^*}$ until the step-size is well adapted. Therefore, $C_T$ and $\tilde{C}_T$ depend on the initial condition, with a logarithmic dependency on the initial multiplicative mismatch.

\subsection{Bounds of the Hitting Time via Drift Conditions}
We are going to use
\emph{drift analysis} that consists in deducing properties of a sequence $\{X_t : t \geq 0 \}$ (adapted to a natural filtration $\{\mathcal{F}_t : t \geq 0\}$) from its drift defined as $\E[X_{t+1} \mid \mathcal{F}_t] - X_t$ \cite{hajek1982hitting}. Drift analysis has been widely used to analyze hitting times of evolutionary algorithms defined on discrete search spaces (mainly on binary search spaces) \cite{he2001drift,he2004study,baritompa1996bounds,mitavskiy2009theoretical,doerr2012multiplicative,doerr2013adaptive}. 
Though they were developed mainly for finite search spaces, the drift theorems can naturally be generalized to continuous domains \cite{lehre2013general,lengler2016drift}. Indeed, Jägersküpper's work \cite{jaegerskuepper2003analysis,jagerskupper2006quadratic,jagerskupper2007algorithmic} is based on the same idea, while the link to the drift analysis was implicit.

Since many drift conditions have been developed for analyzing algorithms on discrete domains, the domain of $X_t$ is often implicitly assumed to be bounded. However, this assumption is violated in our situation, where we will use $X_t = \log\big(\ssf(m_t)\big)$ as the quality measure, which takes values in $\R\cup\{-\infty\}$, and is meant to approach $-\infty$.
We refer to \cite{akimoto2018drift} for additional details.
In general, translating expected progress requires bounding the tail of the progress distribution, as formalized in \cite{hajek1982hitting}.




To control the tails of the drift distribution, we construct a stochastic process $\{Y_t : t \geq 0\}$ iteratively as follows: $Y_0 = X_0$ and
\begin{equation}
   Y_{t+1} =  Y_t + \max \big\{ X_{t+1} - X_t , -A \big\}   \ind{T_\target^X > t} - B   \ind{T_\target^X \leq t}
  \label{eq:truncatedprocess}
\end{equation}
for some $A \geq B > 0$ and $\beta < \beta_0$ with $X_0 = \beta_0$.
It clips $X_{t+1} - X_{t}$ to some constant $-A$ ($A > 0$) from below. We introduce the indicator $\ind{T_\target^X > t}$ for a technical reason. The process disregards progress larger than $A$, and it fixes the progress of the step that hits the target set to $B$.
It is formalized in the following theorem, which is our main mathematical tool to derive an upper bound of the expected first hitting time of \CMAES{} in the form of \eqref{eq:FHT_upper}.
\begin{theorem}
\label{theo:drift-UB-trunc}
Let $\{ X_t : t \geq 0 \}$ be a sequence of real-valued random variables
adapted to a filtration $\{ \F_t : t \geq 0\}$ with $X_0 = \target_0 \in \R$.
For $\target < \target_0$, let $T^X_\target = \inf \left\{ t : X_t \leq \target \right\}$
be the first hitting time of the set $(-\infty, \target]$.
Define a stochastic process $\{Y_t : t \geq 0\}$ iteratively through \eqref{eq:truncatedprocess} with $Y_0 = X_0$ for some $A \geq B > 0$, and let $T^{Y}_\target = \inf \left\{ t :  Y_t \leq \target \right\}$ be the first hitting time of the set $(-\infty, \target]$.
If $ Y_t$ is integrable, i.e.\ $\E\left[\big| Y_t\big|\right] < \infty$, and
\begin{equation}\label{eq:drift-truncated}
\E \left[ \max \left\{X_{t+1} - X_t, -A \right\}  \ind{T_\target^X > t} \,\big|\, \F_t \right] \leq - B  \ind{T_\target^X > t}
	\enspace,
\end{equation}
then the expectation of $T^X_\target$ satisfies
\begin{equation}\label{Bound-HittingTime}
	\E \left[ T^X_\target \right] \leq \E \left[ T^{ Y}_\target \right] \leq \frac{A + \beta_0 - \beta}{B}
	\enspace.
\end{equation}
\end{theorem}
\begin{proof}[Proof of \Cref{theo:drift-UB-trunc}]

  \correct{%
We consider the stopped process $\bar{X}_{t} = X_{\min\{t, T^{X}_\target\}}$. We have $X_t \leq \bar{X}_t$ for $t \leq T^{X}_{\target}$ and $\bar{X}_{t} \leq Y_{\min\{t, T^{X}_\target\}}$ for all $t \geq 0$. Therefore, we have $T^{X}_\target = T^{\bar{X}}_\target \leq T^Y_\target$. Let $\bar{Y}_{t} =  Y_{\min\{t, T^{ Y}_\target\}}$. By construction it holds
$ Y_t \leq \bar{Y}_t$ for $t \leq T^Y_{\target}$ and $T^{Y}_\target = T^{\bar{Y}}_\target$. Hence, $T^{X}_\target \leq T^{Y}_\target \leq T^{\bar{Y}}_\target$.
}{%
  We consider the stopped process $\bar{Y}_{t} =  Y_{\min\{t, T^{ Y}_\target\}}$.
  Then, we have $Y_t = \bar{Y}_t$ for $t \leq T^Y_{\target}$ and $\bar{Y}_t \geq Y_{T^{ Y}_\target}$ for all $t$.
}%

We will prove that
\begin{equation}\label{eq:super-martingale}
	\E[\bar{Y}_{t+1} \mid \F_t ] \leq \bar{Y}_t - B   \ind{T^{ Y}_\target > t} \enspace.
\end{equation}
We start from
\begin{equation}\label{eq:add-two}
	\E[\bar{Y}_{t+1} \mid \F_t ] = \E[ \bar{Y}_{t+1} \ind{ T^{ Y}_\target \leq t } \mid \F_t ] + \E[ \bar{Y}_{t+1} \ind{ T^{Y}_\target > t } \mid \F_t ]
\end{equation}
and bound the different terms:
\begin{equation}\label{eq:top}
	\E[ \bar{Y}_{t+1} \ind{ T^{ Y}_\target \leq t } \mid \F_t ]  = \E[ \bar{Y}_{t} \ind{ T^{ Y}_\target \leq t } \mid \F_t ]  = \bar{Y}_{t} \ind{ T^{ Y}_\target \leq t } \enspace,
\end{equation}
where we have used that $1_{\{T_\target^X > t\}}$, $Y_t$, $1_{\{T_\target^Y > t\}}$, and $\bar{Y}_t$ are all $\F_t$-measurable.
Also
\begin{multline}\label{eq:LEEXtt}
	\E[ \bar{Y}_{t+1} \ind{ T^{ Y}_\target > t } \mid \F_t ]
	= \E[ Y_{t+1} \mid \F_t] \ind{ T^{ Y}_\target > t }
	\\
	\leq (Y_t - B \ind{ T^{X}_\target > t } - B \ind{ T^{X}_\target \leq t }) \ind{ T^{ Y}_\target > t }
	= (\bar{Y}_t - B) \ind{ T^{ Y}_\target > t }
	\enspace,
\end{multline}
where we have used condition \cref{eq:drift-truncated}. Hence, by injecting \cref{eq:top} and \cref{eq:LEEXtt} into \cref{eq:add-two}, we obtain \cref{eq:super-martingale}.
From \cref{eq:super-martingale}, by taking the expectation we deduce
$\E[ \bar{Y}_{t+1}  ] \leq \E[ \bar{Y}_{t}] - B   \Pr[T_{\target}^{Y} > t]$.
Following the same approach as \cite[Theorem~1]{lengler2016drift}, since
$T^{ Y}_\target$ is a random variable taking values in $\mathbb{N}$, it can be
rewritten as $\E[T^{ Y}_\target] = \sum_{t=0}^{+\infty} \Pr[T^{ Y}_\target > t]$
and thus it holds
\begin{equation*}
	B   \E\left[T^{ Y}_\target\right]
	\stackrel{\tilde{t} \to \infty}{\longleftarrow} \sum_{t=0}^{\tilde{t}} B   \Pr\left[T^{ Y}_\target > t\right]
	\leq \sum_{t=0}^{\tilde{t}} \Big( \E[\bar{Y}_t] - \E[\bar{Y}_{t+1}] \Big)
	= \E[\bar{Y}_0] - \E[\bar{Y}_{\tilde{t}}]
	\enspace.
      \end{equation*}
      Since $Y_{t+1} \geq Y_t - A$, we have $Y_{T^{ Y}_\target} \geq \target - A$.
      Given that $\bar{Y}_{\tilde t} \geq Y_{T^{ Y}_\target}$, we deduce that $\E[ \bar{Y}_{\tilde{t}} ] \geq \target - A$
for all $\tilde{t}$. With $\E[\bar{Y}_0] = \beta_0$, we have
$$
\E \left[ T^{ Y}_\target \right] \leq (A / B) + B^{-1}  (\beta_0 - \target) \enspace.
$$
\correct{Since $\E[T^{X}_\target] \leq \E[T^{ Y}_\target]$, this completes the proof.}{%
Because $\target < X_{t} \leq Y_{t}$ for $0 \leq t < T^{X}_\target$, we have $T^{X}_\target \leq T^{Y}_\target$,
implying $\E[T^{X}_\target] \leq \E[T^{ Y}_\target]$. This completes the proof.
}%
\end{proof}
This theorem can be intuitively understood: we assume for the sake of simplicity a process $X_t$ such that $X_{t+1}  \geq X_t -  A$. Then \eqref{eq:drift-truncated} states that the process progresses in expectation by at least $-B$. The theorem concludes that the expected time needed to reach a value smaller than $\beta$ when started in $\beta_0$ equals to $(\beta_0 - \beta)/B$ (what we would get for a deterministic algorithm) plus $A/B$. This last term is due to the stochastic nature of the algorithm. It is minimized if $A$ is as close as possible to $B$, which corresponds to a highly concentrated process.



Jägersküpper \cite[Theorem~2]{jagerskupper2006quadratic} established a general lower bound of the expected first hitting time of the (1+1)-ES. We borrow the same idea to prove the following general theorem for a lower bound of the expected first hitting time, which generalizes \cite[Lemma~12]{jagerskupper2007algorithmic}.
See Theorem~2.3 in \cite{akimoto2018drift} for its proof.
%
\begin{theorem}
\label{theorem:drift-anne-jens}
Let $\{ X_t : t \geq 0 \}$ be a sequence of real-valued random variables
adapted to a filtration $\{ \F_t : t \geq 0\}$ and integrable such that
$X_0 = \target_0$, $X_{t+1} \leq X_t$, and $\E[ X_{t+1} \mid \F_t ] - X_t \geq -C$
for $C > 0$. For $\target < \target_0$ we define
$T^X_\target = \min \left\{ t : X_t \leq \target \right\}$.
Then the expected hitting time is lower bounded by
$\E \left[ T_\target^X \right] \geq - (1/2) + (4 C)^{-1}   (\target_0 - \target)$.
\end{theorem}

\section{Main Result: Expected First Hitting Time Bound}\label{sec:driftanalysis}


\subsection{Mathematical Modeling of the Algorithm}

In the sequel, we will analyze the process $\{ \theta_t : t \geq 0 \}$ where $\theta_t = (m_t,\sigma_t,\Sigma_t) \in \R^n \times \R_> \times \mathcal{S}_\kappa$ generated by the {\CMAES} algorithm. We assume from now on that the optimized objective function $f$ is measurable with respect to the Borel $\sigma$-algebra. We equip the state-space $\mathcal{X} = \R^n \times \R_> \times \mathcal{S}_\kappa$ with its Borel $\sigma$-algebra denoted $\mathcal{B}(\mathcal{X})$.

\subsection{Preliminaries}\label{sec:prelimin}
We present two preliminary results.  In Assumption \ref{asm:radius}, we assume that for $m \in \Xab$, we can include a ball of radius $C_{\ell} \ssf(m)$ into the sublevel set $S_0(m)$ and embed  $S_0(m)$ into a ball of radius $C_u \ssf(m)$. This allows us to upper bound and lower bound the probability of success for all $m \in \Xab$, for all $\Sigma \in \mathcal{S}_\kappa$, by the probability to sample inside of balls of radius $C_u \ssf(m)$ and $C_\ell \ssf(m)$ with appropriate center. From this we can upper-bound $p^\mathrm{upper}_{(a, b]}(\bar\sigma)$ by a function of $\bar\sigma$. Similarly we can lower-bound $p^\mathrm{lower}_{(a, b]}(\bar\sigma)$ by a function of $\bar\sigma$.
The corresponding proof is given in \Cref{apdx:pro:pupper}.
\begin{proposition}\label{pro:pupper}
Suppose that $f$ satisfies \ref{asm:radius}. Consider the lower and upper success probabilities $p^\mathrm{upper}_{(a, b]}$ and $p^\mathrm{lower}_{(a, b]}$ defined in Definition~\ref{def:lusp}, then
\begin{align}
p^\mathrm{upper}_{(a, b]}(\bar\sigma) &\leq \kappa^{d/2} \Phi\left(\bar{\Ball}\left(0, \frac{C_u}{\bar\sigma \kappa^{1/2}}\right); 0, \id \right) \label{eq:pupper}\\
p^\mathrm{lower}_{(a, b]}(\bar\sigma) &\geq \kappa^{-d/2} \Phi\left(\bar{\Ball}\left(\frac{(2 C_u - C_{\ell}) \kappa^{1/2}}{\bar\sigma} e_1, \frac{C_{\ell} \kappa^{1/2}}{\bar\sigma}\right); 0, \id \right) \label{eq:plower}\enspace,
\end{align}
where $e_1 = (1, 0, \dots, 0)$.
\end{proposition}

\newcommand{\pupper}{p^\mathrm{upper}_{m_0}(\bar\sigma)}
\newcommand{\plower}{p^\mathrm{lower}_{m_0}(\bar\sigma)}
We use the previous proposition to establish the next lemma that guarantees the existence of a finite range of normalized step-size that leads to the success probability into some range $(p_u,p_\ell)$ independent of $m$ and $\Sigma$, and provides a lower bound on the success probability with rate $r$ when the normalized step-size is in the above range. Its proof is provided in \Cref{apdx:pro:exist_l_u}.
\begin{lemma}\label{pro:exist_l_u}
We assume that $f$ satisfies \ref{asm:radius} and \ref{asm:levelset} for some $0 \leq a < b \leq \infty$. Then, for any $p_u$ and $p_{\ell}$ satisfying $0 < p_u < p^\mathrm{target} < p_{\ell} < p^\mathrm{limit}$, the constants
	\begin{align*}
		\bar\sigma_{\ell}
			&= \sup \left\{	\bar\sigma > 0 : p^\mathrm{lower}_{(a, b]}(\bar\sigma) \geq p_{\ell} \right\}
			&\text{and}&&
		\bar\sigma_u
			&= \inf \left\{ \bar\sigma > 0 : p^\mathrm{upper}_{(a, b]}(\bar\sigma) \leq p_u \right\}
	\end{align*}
	exist as positive finite values.
	Let $\ell \leq \bar\sigma_{\ell}$ and $u \geq \bar\sigma_u$ such that $u / \ell \geq \aup / \adown$. Then, for $r \in [0,1]$, $p^*_r$ defined as
	\begin{equation}
	p^*_r  :=  \inf_{ \ell \leq \bar\sigma \leq u }  \inf_{m \in \mathcal{X}_a^b} \inf_{\Sigma \in \mathcal{S}_\kappa} p^\mathrm{succ}_{r}(\bar\sigma; m, \Sigma)
	\end{equation}
	is lower bounded by a positive number defined by
	\begin{equation}
          \min_{\ell \leq \bar\sigma \leq u} \kappa^{-d/2} \Phi\left(\Ball\left( \left( \frac{(2 C_u - (1 - r)C_{\ell}) \kappa^{1/2} }{\bar\sigma}\right) e_1, \frac{(1-r) C_{\ell} \kappa^{1/2}}{\bar\sigma} \right); 0, \id \right)
		\enspace.
	\label{eq:p_star}
	\end{equation}
      \end{lemma}



\subsection{Potential Function}\label{sub:potential}

\Cref{pro:exist_l_u} divides the domain of the normalized step-size into three disjoint subsets: $\bar\sigma \in (0, \ell)$ is a too small normalized step-size situation where we have $p^\mathrm{succ}_{0}(\bar\sigma; m, \Sigma) \geq p_{\ell}$ for all $m \in \mathcal{X}_a^b$ and $\Sigma \in \mathcal{S}_\kappa$; $\bar\sigma \in (u, \infty)$ is a too large normalized step-size situation where we have $p^\mathrm{succ}_{0}(\bar\sigma; m, \Sigma) \leq p_u$ for all $m \in \mathcal{X}_a^b$ and $\Sigma \in \mathcal{S}_\kappa$; and $\barsigma \in [\ell, u]$ is a reasonable normalized step-size situation where the success probability with rate $r$ is lower bounded by \cref{eq:p_star}. Since $p_\mathrm{target} \in [p_u, p_\ell]$, the normalized step-size is supposed to be maintained in the reasonable range.

Our potential function is defined as follows. In light of \Cref{pro:exist_l_u}, we can take $\ell \leq \bar\sigma_{\ell}$ and $u \geq \bar\sigma_u$ such that $u / \ell \geq \aup / \adown$. With some constant $v>0$, we define our potential function as
\begin{equation}
  V(\theta) = \log(\ssf(m)) + \max\left\{	0,
    v\log\left[\frac{\aup   \ell   \ssf(m)}{\sigma}\right],
    v\log\left[\frac{\sigma}{\adown   u   \ssf(m)}\right]
  \right\}
  \enspace. \label{eq:potentialfunction}
\end{equation}

The rationale behind the second term on the right-hand side (RHS) is as follows. The second and third terms inside $\max$ are positive only if the normalized step-size $\bar\sigma = \sigma / \ssf(m)$ is smaller than $\ell \aup$ and greater than $u \adown$, respectively. The potential value is $\log\ssf(m)$ if the normalized step-size is in $[\ell \aup, u\adown]$ and it is penalized if the normalized step-size is too small or too large. We need this penalization for the following reason. If the normalized step-size is too small, the success probability is close to $1/2$ for non-critical points, assuming $f= g \circ h$ where $h$ is a continuously differentiable function, but the progress per step is very small because the step-size directly controls the progress for instance measured as $\| m_{t+1} - m_t \| = \sigma_t \|\mathcal{N} (0,\Sigma_t)\|1_{\{f(m_{t+1})\leq f(m_t)\}} $. If the normalized step-size is too large, the success probability is close to zero and produces no progress with high probability. If we would use $\log\ssf(m)$ as a potential function instead of $V(\theta)$ then the progress is arbitrarily small in such situations, which prevents the application of drift arguments. The above potential function penalizes such situations, and guarantees a certain progress in the penalized quantity since the step-size will be increased or decreased, respectively, with high probability, leading to a certain decrease of $V(\theta)$.
We illustrate in \Cref{fig:on_quad} that $\log(\ssf(m))$ cannot work alone as a potential function while $V(\theta)$ does: when we start from a too small or too large step-size, $\log(\ssf(m))$ looks constant (doted line in green and blue). Only when the step-size is started at $1$, we see progress in $\log(\ssf(m))$. Also, the step size can always get arbitrarily worse, with a very small probability, which forces us to handle the case of badly adapted step size properly. Yet the simulation of $V(\theta)$ shows that in all three situations (small, large and well adapted step-sizes compared to the distance to the optimum), we observe a geometric decrease of $V(\theta)$.


\subsection{Upper Bound of the First Hitting Time}\label{sub:upper_bound}
We are now ready to establish that the potential function defined in \eqref{eq:potentialfunction} satisfies a (truncated)-drift condition from \Cref{theo:drift-UB-trunc}. This will in turn imply an upper bound on the expected hitting time of $\ssf(m)$ to $[0,\epsilon]$ provided $a \leq \epsilon$. The proof follows the same line of argumentation as the proof of \cite[Proposition~4.2]{akimoto2018drift}, which was restricted to the case of spherical functions. It was generalized under similar assumptions as in this paper, but for a fixed covariance matrix equal to the identity, in \cite[Proposition~6]{morinaga2019generalized}. The detailed proof is given in \Cref{apdx:pro:drift-UB-trunc}.
\begin{proposition}\label{pro:drift-UB-trunc}
  Consider the \CMAES{} described in \Cref{algo} with state $\theta_t = (m_t, \sigma_t, \Sigma_t)$. Assume that the minimized objective function $f$ satisfies \ref{asm:radius} and \ref{asm:levelset} for some $0 \leq a < b \leq \infty$. Let $p_{u}$ and $p_{\ell}$ be constants satisfying $0 < p_u < \ptarg < p_{\ell} < p^\mathrm{limit}$ and $p_{\ell} + p_u = 2 \ptarg$. Then, there exists $\ell \leq \bar\sigma_{\ell}$ and $u \geq \bar\sigma_u$ such that $u / \ell \geq \aup / \adown$, where $\bar\sigma_{\ell}$ and $\bar\sigma_u$ are defined in \Cref{pro:exist_l_u}. For any $A > 0$, taking $v$ satisfying $0 < v < \min\left\{1,\ \frac{A}{\log(1 /\adown)},\ \frac{A}{\log(\aup)} \right\}$, and the potential function~\cref{eq:potentialfunction}, we have
  \begin{equation}\label{eq:Bdef}
    \E\left[\max\{V(\theta_{t+1})-V(\theta_t), -A\}   \ind{ m_t \in \mathcal{X}_{a}^{b} } \mid \mathcal{F}_t \right] \leq -B   \ind{ m_t \in \mathcal{X}_{a}^{b} }
  \end{equation}
  where
  \begin{equation}
    B = \min\left\{
      A   p^*_r - v   \log\left(\frac{\aup}{\adown}\right),\
      v  \frac{p_{\ell} - p_u}{2}   \log\left(\frac{\aup}{\adown}\right)
    \right\} \enspace,
    \label{eq:B}
  \end{equation}
  and
   $ p^*_r = \inf_{\bar\sigma \in [\ell, u]}
    \inf_{m \in \mathcal{X}_{a}^{b}} \inf_{\Sigma \in \mathcal{S}_\kappa} p^{\mathrm{succ}}_{r}(\bar\sigma; m,\Sigma)
    \ \text{with} \ r = 1-\exp\left(-\frac{A}{1-v}\right).
  $
  Moreover, for any $A > 0$ there exists $v$ such that $B < A$ is positive.
\end{proposition}

We apply \Cref{theo:drift-UB-trunc} along with \Cref{pro:drift-UB-trunc} to derive the expected first hitting time bound. To do so, we need to confirm that it satisfies the prerequisite of the theorem: integrability of the process $\{Y_t : t \geq 0\}$ defined in \cref{eq:truncatedprocess} with $X_t = V(\theta_t)$.

\begin{lemma}\label{lem:int}
	Let $\{\theta_t : t \geq 0\}$ be the sequence of parameters $\theta_t = (m_t, \sigma_t, \Sigma_t)$ defined by the \CMAES{} with the initial condition $\theta_0 = (m_0, \sigma_0, \Sigma_0)$ optimizing a measurable function $f$. Set $X_t = V(\theta_t)$ as defined in \cref{eq:potentialfunction} and define the process $Y_t$ as defined in \Cref{theo:drift-UB-trunc}. Then, for any $A > 0$, $\{Y_t : t \geq 0\}$ is integrable, i.e., $\E[\abs{Y_t}] < \infty$ for each $t$.
\end{lemma}
\begin{proof}[Proof of \Cref{lem:int}]
The drift
$
Y_{t+1} =  Y_t + \max \big\{ V(\theta_{t+1}) - V(\theta_t) , -A \big\}   \ind{T_\target^X > t} - B   \ind{T_\target^X \leq t}
$
is by construction bounded by $-A$ from below. It is also bounded by a constant from above. Indeed, from the proof of \Cref{pro:drift-UB-trunc}, it is easy to find the upper bound, say $C$, of the truncated one-step change, $\Delta_t$ in the proof of \Cref{pro:drift-UB-trunc}, without using \ref{asm:radius} and \ref{asm:levelset}. Let $D = \max\{A, C\}$. Then, by recursion, $\abs{V(\theta_t)} \leq \abs{V(\theta_0)} + \abs{V(\theta_t) - V(\theta_0)} \leq \abs{Y_0} + D   t$. Hence $\E[\abs{Y_t}] \leq \abs{Y_0} + D   t < \infty$ for all $t$.
\end{proof}

Finally, we derive the expected first hitting time of $\log \ssf(m_t)$.
\begin{theorem}\label{theo:FHT-UB}
  Consider the same situation as described in \Cref{pro:drift-UB-trunc}. Let $T_{\epsilon} = \min\{t : \ssf(m_t) \leq \epsilon \}$ be the first hitting time of $\ssf(m_t)$ to $[0, \epsilon]$. Choose $a \leq \epsilon < \ssf(m_t) \leq b$, where $a$ and $b$ appear in \Cref{def:lusp}. 
  If $m_0 \in \Xab$, the first hitting time is upper bounded by	$\mathbb{E}[ T_{\epsilon} ] \leq \big( V(\theta_0) - \log(\epsilon) + A \big) / B$ for $A > B > 0$ described in \Cref{pro:drift-UB-trunc},
	where $V(\theta)$ is the potential function defined in \cref{eq:potentialfunction}.
	Equivalently, we have $\E[T_{\epsilon}] \leq C_T + C_R^{-1}   \log(\ssf(m_0) / \epsilon)$, where
	\begin{align*}
	C_T &= \frac{A}{B} + \frac{v}{B} \max\left\{0,
								\log\left(\frac{\aup   \ell   \ssf(m_0)}{\sigma_0}\right),
								\log\left(\frac{\sigma_0}{\adown   u   \ssf(m_0)}\right)
						\right\}\enspace, &
	C_R &= B \enspace.
	\end{align*}
Moreover, the above result yields an upper bound of the expected first hitting time of $\norm{m_t - x^*}$ to $[0, 2 C_u \epsilon]$.
\end{theorem}
\begin{proof}
\Cref{theo:drift-UB-trunc} with \Cref{pro:drift-UB-trunc} and \Cref{lem:int} together bounds the expected first hitting time of $V(\theta_t)$ to $(-\infty, \log(\epsilon)]$ by $\big( V(\theta_0) - \log(\epsilon) + A \big) / B$. Since $\log\ssf(m_t) \leq V(\theta_t)$, $T_\epsilon$ is bounded by the first hitting time of $V(\theta_t)$ to $(-\infty, \log(\epsilon)]$. The inequality is preserved if we take the expectation. The last claim is trivial from the inequality $\norm{x - x^*} \leq 2 C_u \ssf(x)$, which holds under \ref{asm:radius}.
\end{proof}

\Cref{theo:FHT-UB} shows an upper bound on the expected hitting time of the \CMAES{} with success-based step-size
adaptation for linear convergence towards the global optimum $x^*$ on functions
satisfying \ref{asm:radius} and \ref{asm:levelset} with $a = 0$.
Moreover, for $b = \infty$, this bound holds from all initial search points $m_0$.
If $a > 0$, the bound in Theorem~\ref{theo:FHT-UB} does not translate into linear convergence,
but we still obtain an upper bound on the expected first hitting time of
the target accuracy $\epsilon \geq a$. This is useful for understanding the
behavior of \CMAES{} on multimodal functions, and on functions with
degenerated Hessian matrix at the optimum.

\subsection{Lower Bound of the First Hitting Time}\label{sub:lower_bound}

We derive a general lower bound of the expected first hitting time of $\norm{m_t - x^*}$ to $[0, \epsilon]$. The following results hold for an arbitrary measurable function $f$ and for a \CMAES{} with an arbitrary $\sigma$-control mechanism. The following lemma provides the lower bound of the expected one-step progress measured by the logarithm of the distance to the optimum.

\begin{lemma}\label{lem:sphere-progress-bound}
We consider the process $\{\theta_t : t \geq 0\}$ generated by a \CMAES{} algorithm with an arbitrary step-size adaptation mechanism and an arbitrary covariance matrix update optimizing an arbitrary measurable function $f$. We assume $d \geq 2$ and $\kappa_t = \Cond(\Sigma_t) \leq \kappa$. We consider the natural filtration $\F_t$.
Then, the expected single-step progress is lower-bounded by
\begin{equation}
\E[\min(\log(\norm{m_{t+1} - x^*} / \norm{m_t - x^*}), 0) \mid \F_{t}]  \geq - \kappa_t^{\frac{d}{2}} / d \enspace.
\end{equation}
\end{lemma}
\begin{proof}[Proof of \Cref{lem:sphere-progress-bound}]
Note first that
$\log( \norm{m_{t+1} - x^*} / \norm{m_{t} - x^*} ) = \log( \norm{x_t - x^*} / \norm{m_{t} - x^*} ) \indlr{f(x_t) \leq f(m_t)}$.
This value can be positive since $f(x_t) \leq f(m_t)$ does not imply $\norm{x_t - x^*} \leq \norm{m_t - x^*}$ in general. Clipping the positive part to zero, we obtain a lower bound, which is the RHS of the above equality times the indicator $\indlr{\norm{x_t - x^*} \leq \norm{m_t - x^*}}$. Since the quantity is non-positive, dropping the indicator $\indlr{f(x_t) \leq f(m_t)}$ only decreases the lower bound. Hence, we have
$\min(\log( \norm{m_{t+1} - x^*} / \norm{m_{t} - x^*} ), 0) \geq \log( \norm{x_t - x^*} / \norm{m_{t} - x^*} )\indlr{\norm{x_t - x^*} \leq \norm{m_t - x^*}}$.
Then,
\begin{multline*}
\E[\min(\log(\norm{m_{t+1} - x^*}) - \log(\norm{m_t - x^*}), 0) \mid \F_{t}]  \\
\geq \E[ \log( \norm{x_t - x^*} / \norm{m_{t} - x^*} )\indlr{\norm{x_t - x^*} \leq \norm{m_t - x^*}} \mid \F_{t}]
\enspace.
\end{multline*}

We rewrite the lower bound of the drift. The RHS of the above inequality is the integral of $\log(\norm{x - x^*} / \norm{m_t - x^*})$ in the integral domain $\bar\Ball(x^*, \norm{m_t - x^*})$ under the probability measure $\Phi\left(  ; m_t, \sigma_t^2 \Sigma_t \right)$. Performing a variable change (through rotation and scaling) so that $m_t - x^*$ becomes $e_1 = (1, 0, \cdots, 0)$
and letting $\tilde\sigma_t = \sigma_t / \norm{m_t - x^*}$, we can further rewrite it as the integral of $\log(\norm{x})$ in $\bar\Ball(0, 1)$ under $\Phi\left(  ; e_1, \tilde\sigma_t^2 \Sigma_t \right)$. With $\kappa_t = \Cond(\Sigma_t)$, we have
$\varphi\left(  ; e_1, \tilde\sigma_t^2 \Sigma_t \right) \leq \kappa_t^{d/2} \varphi\left(  ; e_1, \kappa_t \tilde\sigma_t^2 \id \right)$, see \Cref{lemma:phi}. 
Altogether, we obtain the lower bound
$\E[ \log( \norm{x_t - x^*} / \norm{m_{t} - x^*} )\indlr{\norm{x_t - x^*} \leq \norm{m_t - x^*}} \mid \F_{t}]
\geq \kappa_t^{d/2} \int_{\bar\Ball(0, 1)} \log(\norm{x}) \varphi\left(  ; e_1, \kappa_t \tilde\sigma_t^2 \id \right) \mathrm{d}x$.
The RHS is equivalent to $-\kappa_t^{d/2}$ times the single step progress of the (1+1)-ES on the spherical function at $m_t = e_1$ and $\sigma = \sqrt{\kappa} \tilde\sigma_t$%
, which is proven in the proof of Lemma~4.4 of \cite{akimoto2018drift} to be lower bounded by $1/d$ for $d \geq 2$. This completes the proof.
\end{proof}

The following theorem proves that the expected first hitting time of \CMAES{} is $\Omega(\log(\norm{m_0 - x^*} / \epsilon))$ for any measurable function $f$, implying that it can not converge faster than linearly. In case of $\kappa = 1$ the lower runtime bound becomes $\Omega(d   (\log(\norm{m_0 - x^*} / \epsilon)) )$, meaning that the runtime scales linearly with respect to $d$. The proof is a direct application of  \Cref{lem:sphere-progress-bound} to \Cref{theorem:drift-anne-jens}.

\begin{theorem}\label{thm:eslower}
We consider the process $\{\theta_t : t \geq 0\}$ generated by a \CMAES{} described in \Cref{algo} and assume that $f$ is a measurable function with $d \geq 2$. Let $T_\epsilon = \inf \{ t : \| m_t - x^*\| \leq \epsilon \}$ be the first hitting time of $[0, \epsilon]$ by $\|m_t - x^*\|$. Then, the expected first hitting time is lower bounded by
$\E[T_{\epsilon}] \geq -( 1/ 2) + \frac{d}{4 \kappa^{d/2}} \log (\norm{m_0 - x^*} / \epsilon)$.
The bound holds for arbitrary step-size adaptation mechanisms. If \ref{asm:radius} holds, it gives a lower bound for the expected first hitting time bound of $\ssf(m_t)$ to $[0, 2 C_\ell \epsilon]$.
\end{theorem}
\begin{proof}[Proof of \Cref{thm:eslower}]
Let $X_t = \log\norm{m_t - x^*}$ for $t \geq 0$. Define $Y_{t}$ iteratively as $Y_{0} = X_{0}$ and $Y_{t+1} = Y_{t} + \min(X_{t+1} - X_{t}, 0)$. Then, it is easy to see that $Y_t \leq X_t$ and $Y_{t + 1} \leq Y_{t}$ for all $t \geq 0$. Note that $\E[Y_{t+1} - Y_{t}\mid \mathcal{F}_t] = \E[\min(X_{t+1} - X_{t}, 0) \mid \mathcal{F}_t] = \E[\min(\log(\norm{m_{t+1} - x^*} / \norm{m_t - x^*}), 0) \mid \F_{t}]$, where the RMS is lower bounded in light of \Cref{lem:sphere-progress-bound}. Then, applying \Cref{theorem:drift-anne-jens}, we obtain the lower bound. The last statement directly follows from $\norm{x - x^*} \leq 2 C_\ell \ssf(x)$ under \ref{asm:radius}.
\end{proof}
\subsection{Almost Sure Linear Convergence}\label{sec:aslinear}

Additionally to the expected first hitting time bound, we can deduce from \Cref{pro:drift-UB-trunc}, almost sure linear convergence as stated in the following proposition.%
\begin{proposition}\label{prop:cr}
  Consider the same situation as described in \Cref{pro:drift-UB-trunc}, where $a = 0$ and $0 < b \leq \infty$.
  Then, for any $m_0 \in \mathcal{X}_{0}^{b}$, $\sigma_0 > 0$ and $\Sigma \in \mathcal{S}_\kappa$, we have
  \begin{equation}\label{eq:as-cv}
  \Pr\left[\limsup_{t\to\infty} \frac1t \log \ssf(m_t) \leq -B \right]
  = \Pr\left[\limsup_{t\to\infty} \frac1t \log \|m_t - x^*\| \leq -B \right]
  = 1\enspace,
\end{equation}
  where $B > 0$ is as defined in \Cref{pro:drift-UB-trunc}.
  Hence almost sure linear convergence holds at a rate $\exp(-C)$ such that $\exp(-C) \leq \exp(-B)$.
\end{proposition}
\begin{proof}[Proof of \Cref{prop:cr}]
  Let $V$ be defined in \eqref{eq:potentialfunction}.
  Let $Y_{0} = V(\theta_0)$ and $Y_{t+1} = Y_t + \max(-A, V(\theta_{t+1}) - V(\theta_{t}))$.
  Define $Z_{t} = Y_{t} - \E_{t-1}[Y_t]$ for $t \geq 0$.
  Then, $\{Z_{t}\}$ is a martingale difference sequence on the filtration $\{\mathcal{F}_t\}$ produced by $\{\theta_t\}$.
  We hence have $\frac1t \log \ssf(m_t) \leq \frac1t V(\theta_t) \leq \frac1t Y_t$, and from \Cref{pro:drift-UB-trunc} we obtain
\begin{align*}
  Y_t
  = \E_{t-1}[Y_t] + Z_{t}
  = Y_{t-1} + \E_{t-1}[Y_t - Y_{t-1}] + Z_{t}
  \leq Y_{t-1} - B  + Z_{t}
  \enspace.
\end{align*}
By repeatedly applying the above inequality and dividing it by $t$, we obtain
$\frac1t Y_t \leq - B  + \frac1t Y_{0} + \frac1t \sum_{i=1}^{t} Z_{i}$,
where $\lim_{t\to\infty} \frac1t Y_{0} = 0$ and $\sum_{i=1}^{t} Z_{i}$ is a martingale sequence.
In light of the strong law of large numbers for martingales \cite{chow1967}, if $\sum_{t=1}^{\infty} \E[Z_t^2] / t^2 < \infty$, we have $\lim_{t \to \infty} \frac1t \sum_{i=1}^{t} Z_{i} = 0$ almost surely.
By the definition of $V(\theta_t)$ and the working mechanism of the \CMAES, we have $V(\theta_{i}) - V(\theta_{i-1}) \leq v \log(\aup / \adown)$. Hence,
  $\E[Z_i^2] = \E[(Y_i - \E_{i-1}[Y_i])^2] = \E[\max(-A, V(\theta_{i}) - V(\theta_{i-1}))^2] \leq \max(A, v \log(\aup / \adown))^2$.
Hence, we have
  $\limsup_{t\to\infty} \frac1t \log \ssf(m_t)
  \leq - B + \lim_{t\to\infty} \frac1t Y_{0} + \lim_{t\to\infty} \frac1t \sum_{i=1}^{t} Z_{i}
  = -B$ almost surely.
Along with $\norm{x - x^*} \leq 2 C_u \ssf(x)$, we obtain \Cref{eq:as-cv}.
\end{proof}

\subsection{Wrap-up of the Results: Global Linear Convergence}\label{sec:func}


As a corollary to the lower-bound from \Cref{thm:eslower}, the upper bound from \Cref{theo:FHT-UB}, Proposition~\ref{prop:cr} stating the almost sure linear convergence and the fact that different assumptions discussed in Section~\ref{sec:assump} imply \ref{asm:radius} and \ref{asm:levelset}, we summarize our linear convergence results in the following theorem.
\begin{theorem}[Global Linear Convergence]
\label{theo:main-result}
We consider the \CMAES{} optimizing an objective function $f$.
Suppose either
\begin{itemize}
\item[(a)]
	$f$ satisfies \ref{asm:radius} and \ref{asm:levelset} for $a=0$, $p^\mathrm{limit} > p^\mathrm{target}$, and $m_0 \in \mathcal{X}_0^b$; or

      \item[(b)]
	$f$ satisfies either \ref{quad-bound} or \ref{asm:pos-hom}, $p^\mathrm{target} < 1/2$, and $m_0 \in \mathbb{R}^d$.
      \end{itemize}
Then, for any $\sigma_0 > 0$ and $\Sigma_0 \in \mathcal{S}_\kappa$, the expected hitting time $\E[T_\epsilon]$ of $\norm{m_t - x^*}$ to $[0, \epsilon]$  is $\Theta\big(\log(\norm{m_0 - x^*} / \epsilon)\big)$ for all $\epsilon > 0$.
Moreover, both $\ssf(m_t)$ and $\norm{m_t - x^*}$ linearly converge almost surely, i.e.
$$
\Pr\left[\limsup_{t\to\infty} \frac1t \log \ssf(m_t) \leq -B \right]
= \Pr\left[\limsup_{t\to\infty} \frac1t \log \norm{m_t-x^*} \leq -B \right] = 1\enspace,
$$
where $B > 0$ is as defined in \Cref{pro:drift-UB-trunc}. The convergence rate $\exp(-C)$ is thus upper-bounded by $\exp(-B)$.
\end{theorem}%

\subsection{Tightness in the Sphere Function Case}\label{sub:sphere}

Now we consider a specific convex quadratic function, namely the sphere function $f(x) = \frac12 \norm{x}^2$ where the spatial suboptimality function equals $\ssf(x) = V_d \norm{x}$.
In \Cref{theo:main-result} we have formulated that the expected hitting time of a ball of radius $\epsilon$ for the \CMAES{} equals $\Theta(\log\|m_0 - x^*\|/\epsilon)$. Yet, this statement does not give information on how the constants hidden in the $\Theta$-notation scale with the dimension.
In particular the convergence rate of the algorithm is upper-bounded by
  $\exp(-B)$
  where $B$ is given in \eqref{eq:B}, see Theorem~\ref{theo:FHT-UB}.
In this section, we estimate precisely the scaling of $B$ in \Cref{pro:drift-UB-trunc} with respect to the dimension and compare it with the general lower bound of the expected first hitting time given in Theorem~\ref{thm:eslower}. We then conclude that the bound is tight with respect to the scaling with $d$ in the case of the sphere function.

Let us assume $\kappa = 1$, that is, we consider the (1+1)-ES without covariance matrix adaptation ($\Sigma = I$). Then, $p^\mathrm{lower}_{(a, b]}(\bar\sigma) = p^\mathrm{upper}_{(a, b]}(\bar\sigma) = p_r^{\mathrm{succ}}(\barsigma; m, \Sigma)$, where the right-most side is independent of $m$ and $\Sigma$ as described in \Cref{lemma:r-proba-success-sphere}. 
This means that the success probability is solely controlled by the normalized step-size $\bar\sigma$.

The following proposition states that the convergence speed is $\Omega(1/d)$, hence the expected first hitting time scales as $O(1/d)$.
The proof is provided in \Cref{apdx:proposition:scaling}.
\begin{proposition}
\label{proposition:scaling}
For $A = 1/d$, $p_\mathrm{target} \in \Theta(1)$ and $\log(\aup / \adown) \in \omega(1/d)$, we have $B \in \Omega(1/d)$.
\end{proposition}

Two conditions on the choice of $\aup$ and $\adown$: $p_\mathrm{target} = \log(1/\adown) / \log(\aup / \adown) \in \Theta(1)$ and $\log(\aup / \adown) \in \omega(1/d)$, are understood as follows. The first condition implies that the target success probability $p_\mathrm{target}$ must be independent of $d$. In the 1/5 success rule, $\aup$ and $\adown$ are set so that $p_\mathrm{target} = 1/5$ independent of $d$. The second condition implies that the factors of the step-size increase and decrease must be $\log(\aup) \in \omega(1/d)$ and $\log(1/\adown) \in \omega(1/d)$. Note that on the sphere function the normalized step-size $\barsigma \propto \sigma / \norm{m - x^*}$ is kept around a constant during the search. It implies that the convergence speed of $\norm{m - x^*}$ and $\sigma$ must agree. Therefore the speed of the adaptation of the step-size must not be too small to achieve $\Theta(d)$ scaling of the expected first hitting time.

\Cref{proposition:scaling} and \Cref{theo:FHT-UB} imply $\E[T_\epsilon] \in O(d   \log (\norm{m_0} / \epsilon))$ and \Cref{thm:eslower} implies $\E[T_\epsilon] \in \Omega( d   \log (\norm{m_0} /\epsilon) )$. They yield $\E[T_\epsilon] \in \Theta(d   \log (\norm{m_0} /\epsilon))$. This result shows i) that the runtime of the (1+1)-ES on the sphere function is proportional to $d$ as long as $\log(\aup / \adown) \in \omega(1/d)$, and ii) that from our methodology one can derive a tight bound of the runtime in some cases. The result is formally stated as follows.

\begin{theorem}\label{thm:sphere}
The (1+1)-ES (\Cref{algo}) with $\kappa = 1$ and $p^\mathrm{target} < 1/2$ converges globally and linearly in terms of $\log \norm{m_t - x^*}$ from any starting point $m_0 \in \R^d$, $\sigma_0 > 0$, and $\Sigma_0 = I$ on any function $f(x) = g(\norm{x - x^*})$, where $g$ is a strictly increasing function. Moreover, if $p^\mathrm{target} \in \Theta(1)$ and $\log(\aup / \adown) \in \omega(1/d)$,
the expected first hitting time $T_\epsilon$ of $\log\norm{m_t - x^*}$ to $(-\infty, \log(\epsilon)]$ is $\Theta( d   \log (\norm{m_0} /\epsilon))$ and the almost sure convergence rate is upper-bounded by $\exp(-\Theta(1/d))$.
\end{theorem}

Since the lower bound holds for an arbitrary $\sigma$-adaptation mechanism, the above result not only implies that our upper bound is tight, but it also implies that the success-based $\sigma$-control mechanism achieves the best possible convergence rate except for a constant factor on the spherical function.

\section{Discussion}
\label{sec:discussion}
We have established the almost sure global
linear convergence of the \CMAES{} and also expressed as a bound on the expected
hitting time of an $\epsilon$-neighborhood of the solution. Assumption
\ref{asm:radius} has been the key to obtaining the expected first
hitting time bound of \CMAES{} in the form of \eqref{eq:FHT_upper}. The
convergence results hold on a wide class of functions. It includes
\begin{itemize}
\item[(i)] strongly convex functions with Lipschitz gradient, where
	linear convergence of numerical optimization algorithm is usually
	analyzed,
\item[(ii)] continuously differentiable positively homogenous functions,
	where previous linear convergence results had been introduced, and
\item[(iii)] functions with non-smooth level sets as illustrated in
	Figure~\ref{fig:levelset-balls}.
\end{itemize}
Because the analyzed algorithms are invariant to strictly monotonic
transformations of the objective functions, \emph{all results that hold on $f$
also hold on $g \circ f$ where $g: {\rm Im}(f) \to \mathbb{R}$ is a strictly increasing
transformation, which can thus introduce discontinuities on the objective function}.
In contrast to the previous result establishing the convergence of CMA-ES
\cite{Diouane2015} by adding a step to enforce a sufficient decrease (which works well for direct search methods, but which is unnatural for ESs), we did not need to modify
the adaptation mechanism of the \ES{} to achieve our convergence proofs.
We believe that this is crucial, since it allows our analysis to
reflect the main mechanism that makes the algorithm work well in practice.


\Cref{thm:sphere} proves that we can derive a tight convergence rate with \Cref{pro:drift-UB-trunc} on the sphere function in the case where $\kappa = 1$, i.e., without covariance matrix adaptation. This partially supports the utility of our methodology. However, its derivation relies on the fact that both the level sets of the objective function and the equal-density curves of the sampling distribution are isotropic, and hence does not generalize immediately. Moreover, the lower bound (\Cref{thm:eslower}) seems to be loose even for $\kappa = 1$ on convex quadratic functions, where we empirically observe that the logarithmic convergence rate scales like $\Theta(1 / \Cond(\nabla \nabla f))$, see \Cref{fig:on_quad}, while its dependency on the dimension is tight.

A better lower bound of the expected first hitting time and a handy way to estimate the convergence rate are relevant directions of future work. Further directions of future work are as follows:


Proving linear convergence of \CMAES{} does not reveal the
benefits of \CMAES{} over the \ES{} without covariance matrix
adaptation. The motivation of the introduction of the covariance matrix
is to improve the convergence rate and to broaden the class of functions
on which linear convergence is exhibited. None of them are achieved in
this paper.

On convex quadratic functions, we empirically observe that the covariance matrix approaches a stable distribution that is closely concentrated around the inverse Hessian up to a scalar factor, and the convergence speed on all convex quadratic functions is equal to that on the sphere function (see \Cref{fig:on_quad}). This behavior is not described by our result.

Covariance matrix adaptation is also important for optimizing functions with non-smooth level sets. On continuously differentiable functions, we can always set $\aup$ and $\adown$ so that $p = \frac{\log (1/\adown)}{\log(\aup / \adown)} < p^\mathrm{limit} = 1/2$. This is the rationale behind the 1/5 success rule, where $p = 1/5$. Indeed, $p = 1/5$ is known to approximate the optimal situation on the sphere function where the expected one-step progress is maximized \cite{rechenberg:1973}. Therefore, one does not need to tune these parameters in a problem-specific manner. However, if the objective is not continuously differentiable and levelsets are non-smooth, then $p^\mathrm{limit}$ is in general smaller than $1/2$. For example, it can be as low as $p^\mathrm{limit} = 1/2^d$ on $f(x) = \|x\|_\infty = \max_{i =1,\dots,n}\abs{x_i}$.
Without an appropriate adaptation of the covariance matrix the success probability will be smaller than $p = 1/5$ and one must tune $\aup$ and $\adown$ in order to converge to the optimum, which requires information about $p^\mathrm{limit}$. By adapting the covariance matrix appropriately, the success probability can be increased arbitrary close to $1/2$ (by elongating steps in the direction of the success domain) and $\aup$ and $\adown$ do not require tuning.

To achieve a reasonable convergence rate bound and broaden the class of
functions on which linear convergence is exhibited, one needs to find
another potential function $V$ that may penalize a high condition number
$\Cond(\nabla \nabla f(m_t) \Sigma_t)$ and replace the definitions
of $p^\mathrm{upper}$ and $p^\mathrm{lower}$ accordingly. This point is
left for future work.



\section*{Acknowledgement}
We gratefully acknowledge support by Dag\-stuhl seminar 17191
``Theory of Randomized Search Heuristics''. We would like to thank
Per Kristian Lehre, Carsten Witt, and Johannes Lengler for valuable
discussions and advice on drift theory.
Y.\ A.\ is supported by JSPS KAKENHI Grant Number 19H04179.


\appendix
\let\displaystyle\textstyle
\section{Some Numerical Results}\label{app:NumericalResults}

We present experiments with five algorithms on two convex quadratic
functions. We compare (1+1)-ES, (1+1)-CMA-ES, simplified direction
search \cite{konevcny2014simple}, random pursuit
\cite{stich2013optimization}, and gradientless descent \cite{Golovin2020Gradientless}.

All algorithms were started at the initial search point
$x_0 = \frac{1}{\sqrt{d}} (1, \dots, 1) \in \R^d$. We implemented the
algorithms as follows, with their parameters tuned where necessary:
	The ES always uses the setting $\aup=\exp(4/d)$ and $\adown=\aup^{-1/4}$
	for step size adaptation.
	We set the constant $c$ in the sufficient decrease condition of
	Simplified Direction Search to $\frac{1}{10}$, and we employed the
	standard basis as well as the negatives of these vectors as
	candidate directions. In each iteration we looped over the set of
	directions in random order. Randomizing the order greatly boosted
	performance over a fixed order.
	Random Pursuit was implemented with a golden section line search in
	the range $[-2 \sigma, 2 \sigma]$ with a rather loose target
	precision of $\sigma/2$, where $\sigma$ is either the initial step
	size or the length of the previous step.
	For Gradientless Descent we used the initial step size as the
	maximal step size and defined a target precision of $10^{-10}$. This
	target is reached by the ES in all cases.
%
The experiments are designed to demonstrate several different effects:
(a)
	We perform all experiments in $d=10$ and $d=50$ dimensions to
	investigate dimension-dependent effects.
(b)
	We investigate best-case performance by running the algorithms on
	the spherical function $\|x\|^2$, i.e., on the separable convex quadratic
	function with minimal condition number. The initial step size is set
	to $\sigma_0 = 1$. All algorithms have a budget of $100   d$
	function evaluations.
(c)
	We investigate the dependency of the performance on initial parameter
	settings by repeating the same experiment as above, but with an initial
	step size of $\sigma_0 = \frac{1}{1000}$. All algorithms have a budget
	of $700   d$ function evaluations.
(d)
	We investigate the dependence on problem difficulty by running the
	algorithms on an ellipsoid problem with a moderate condition number
	of $\kappa_f = 100$. The eigenvalues of the Hessian are evenly
	distributed on a log-scale. We use $\sigma_0 = 1$ like in the first
	experiment. All algorithms have a budget of $500   d$ function
	evaluations.
\begin{figure}
\begin{center}
  \includegraphics[width=0.33\hsize]{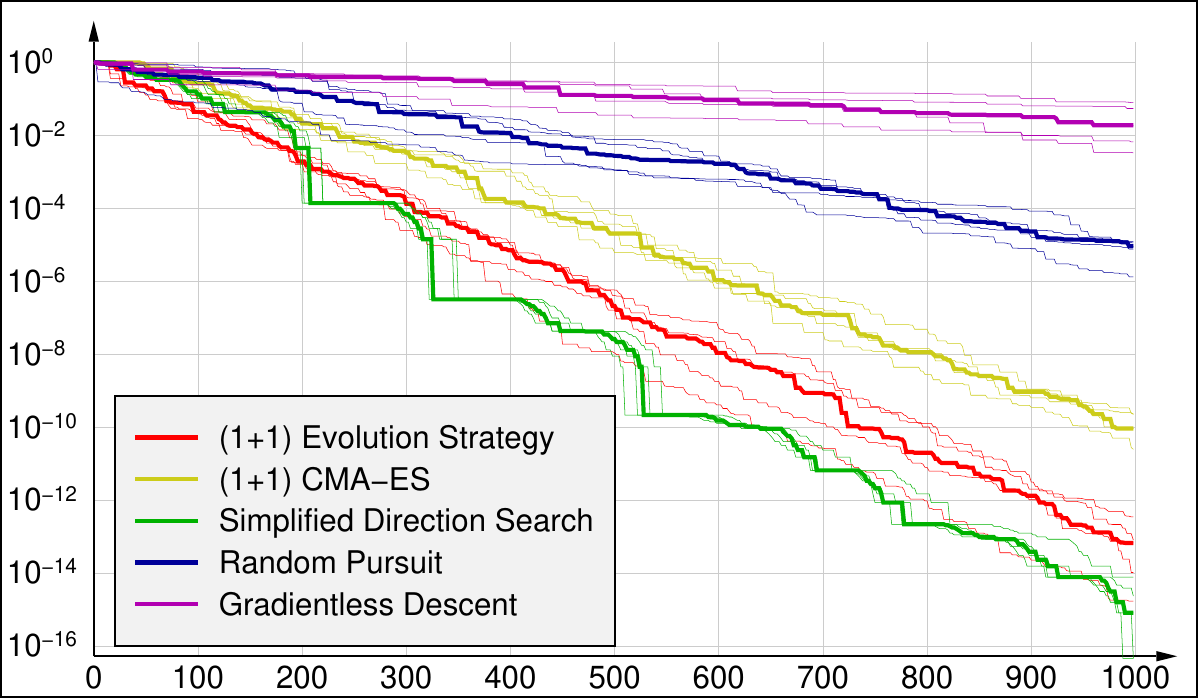}%
  \includegraphics[width=0.33\hsize]{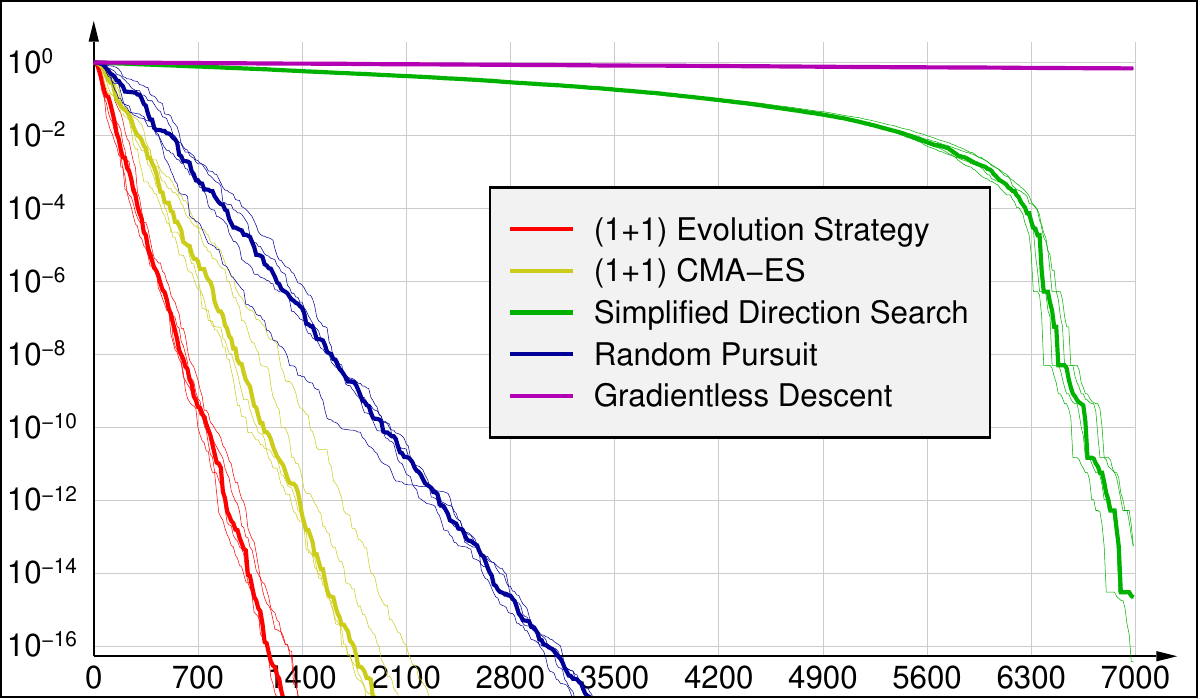}%
  \includegraphics[width=0.33\hsize]{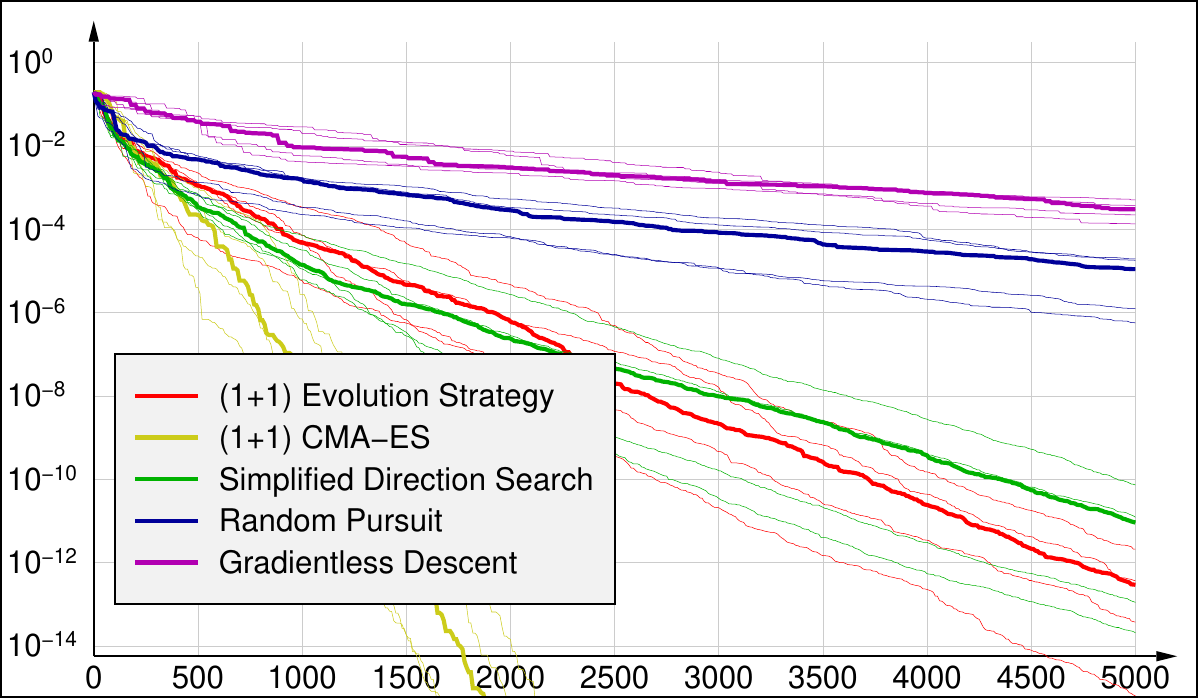}%
  \\[0.5em]
  \includegraphics[width=0.33\hsize]{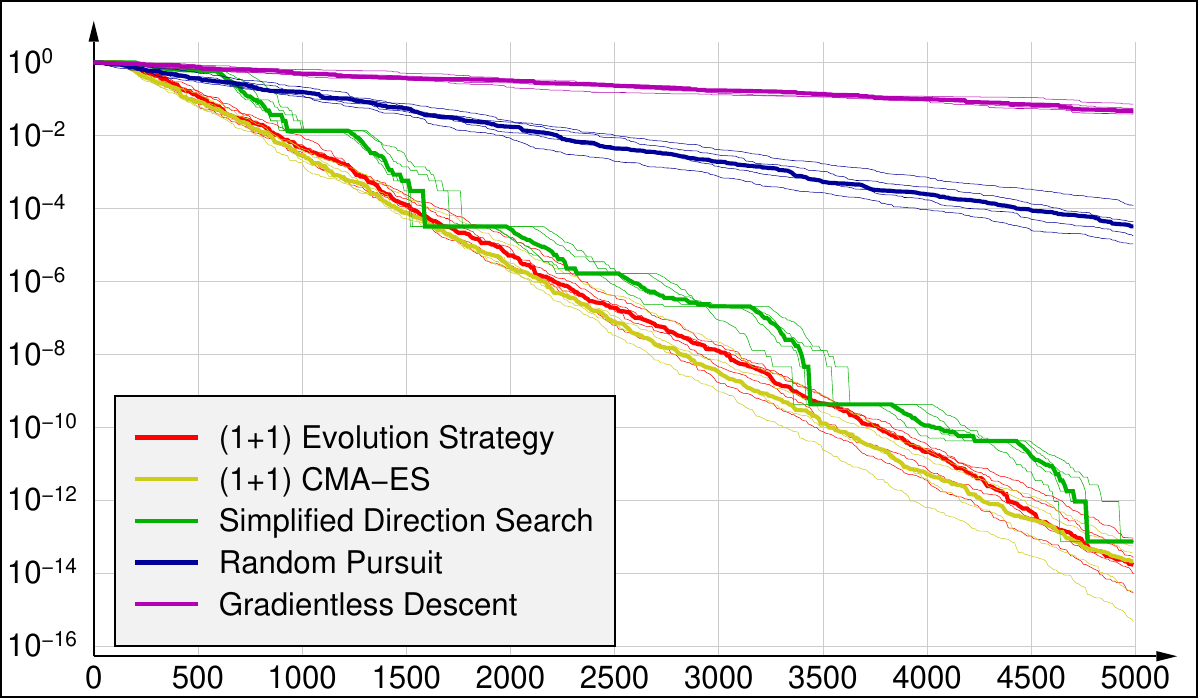}%
  \includegraphics[width=0.33\hsize]{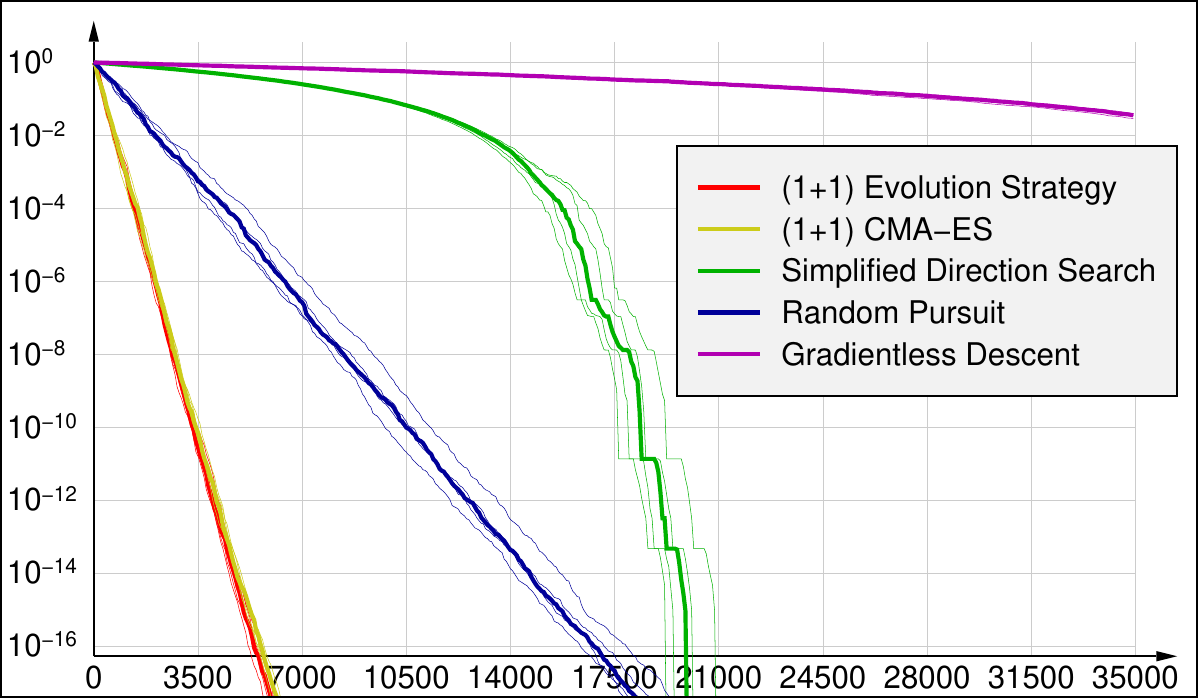}%
  \includegraphics[width=0.33\hsize]{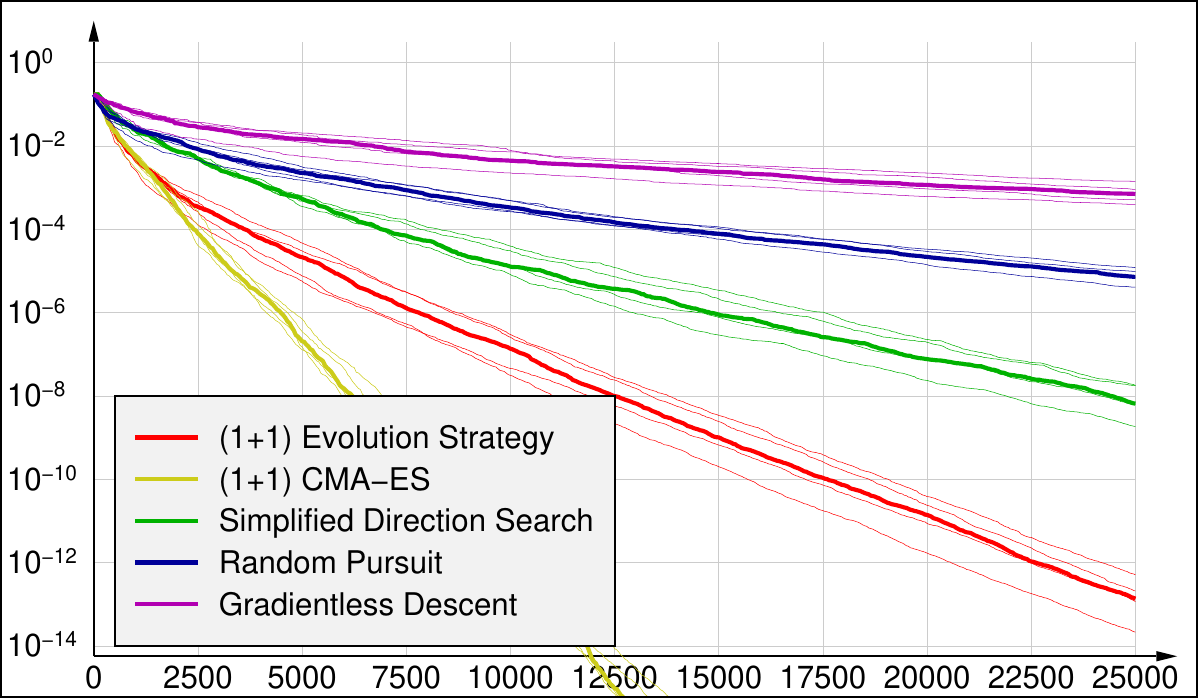}
\caption{
	\label{fig:performance}
	Comparison of (1+1)-ES with and without covariance matrix adaptation
	with three well-analyzed derivative-free optimization algorithms on two convex quadratic functions.
	The left column of plots shows the performance on the sphere function $\|x\|^2$
	in dimensions 10 (top) and 50 (bottom). The middle column shows the same problem,
	but the initial step size is smaller by a factor of $1000$ (and the horizontal
	axis differs), simulating that the distance to the optimum was under-estimated.
	The right column shows the performance on the ellipsoid function (defined in \Cref{fig:on_quad}).
	The plots show the evolution of the best-so-far function value (on a logarithmic
	scale), with five individual runs (thin curves) as well as median performance
	(bold curves).
}
\end{center}
\end{figure}
The experimental results are presented in Figure~\ref{fig:performance}.

\textbf{Interpretation.}
We observe only moderate dimension-dependent effects, besides the
expected linear increase of the runtime.
We see robust performance of the ES, in particular with covariance
matrix adaptation. The second experiment demonstrates the practical
importance of the ability to grow the step size: the ES is essentially
unaffected by wrong initial parameter settings while the gradientless descent and the simplified direct search are (which can be understood directly from the algorithms themselves). This property does not
show up in convergence rates and is therefore often (but not always)
neglected in algorithm design. The last experiment clearly demonstrates
the benefit of variable-metric methods like CMA-ES. It should be noted
that variable metric techniques can be implemented into most existing
algorithms. This is rarely done though, with random pursuit being a
notable exception~\cite{stich2016variable}.

\section{Proofs}\label{appendix}

\subsection{Proof of \Cref{lem:asm1}}\label{apdx:lem:asm1}
Since $\ssf$ is invariant to $g$, without loss of generality we assume $f(x) = h(x) - h(x^*)$ in this proof.
Inequality \cref{eq:quad-bound} implies that $f(y) \leq f(x) \Rightarrow (L_{\ell} / 2) \norm{y - x^*}^2 \leq f(x)$, meaning that $\{y : f(y) \leq f(x) \} \subseteq \bar{\Ball}\Big(x^*, \sqrt{\frac{f(x)}{L_{\ell}/2}}\Big)$.
Since $\ssf(x)$ is the $d$th root of the volume of the left-hand side of the above relation, we find
$\ssf(x) \leq \mu^{\frac1d}\Big(\bar{\Ball}\Big(x^*, \sqrt{\frac{f(x)}{L_{\ell}/2}}\Big)\Big) = V_d\sqrt{\frac{f(x)}{L_{\ell}/2}}$.
Analogously, we obtain $\Ball\Big(x^*, \sqrt{\frac{f(x)}{L_u/2}}\Big) \subseteq \{y : f(y) < f(x) \}$ and $\ssf(x) \geq V_d \sqrt{\frac{f(x)}{L_u/2}}$.
From these inequalities, we obtain $\{ y : f(y) \leq f(x) \} \subseteq \bar{\Ball}\Big(x^*, \sqrt{\frac{L_u}{L_{\ell}}} \frac{\ssf(x)}{V_d}  \Big)$ and $\Ball\Big(x^*,  \sqrt{\frac{L_{\ell}}{L_u}}\frac{\ssf(x)}{V_d} \Big) \subseteq \{ y : f(y) < f(x) \}$.
This implies \ref{asm:radius} for $\mathcal{X}_{0}^{\infty}$.
\ref{asm:levelset} is immediately implied by \Cref{prop:asm2}.
This completes the proof.

\subsection{Proof of \Cref{prop:pos-hom}}\label{apdx:prop:pos-hom}

We first prove that \ref{asm:radius} holds for $a = 0$ and $b = \infty$ with $C_u = \sup \{\norm{x - x^*} : \ssf(x) = 1 \}$ and $C_{\ell} = \inf \{\norm{x - x^*} : \ssf(x) = 1 \}$ and they are finite.

It is easy to see that the spatial suboptimality function $\ssf(x)$ is proportional to $h(x) - h(x^*)$. Let $\ssf(x) = c   (h(x) - h(x^*))$ for some $c > 0$. Then, $\ssf$ is also a homogeneous function. Since it is homogeneous, \ref{asm:radius} reduces to that there are open and closed balls with radius $C_{\ell}$ and $C_u$ satisfying the conditions described in the assumption with \new{$f_\mu(m) = 1$}. Such constants are obtained by $C_u = \sup \{\norm{x - x^*} : \ssf(x) = 1 \}$ and $C_{\ell} = \inf \{\norm{x - x^*} : \ssf(x) = 1 \}$.

Due to the continuity of $f$ there exists an open ball $\mathcal{B}$ around $x^*$ such that $h(x) < h(x^*) + 1/c$ for all $x \in \mathcal{B}$. Then, it holds that $\ssf(x) < 1$ for all $x \in \mathcal{B}$. It implies that $C_{\ell}$ is no smaller than the radius of $\mathcal{B}$, which is positive. Hence, $C_{\ell} > 0$.

We show the finiteness of $C_u$ by a contradiction argument. Suppose $C_u = \infty$. Then, there is a direction $v$ such that $\ssf(x^* + M  v) \leq 1$ with an arbitrarily large $M > 0$. Since $\ssf$ is homogeneous, we have $\ssf(x^* + v) \leq 1 / M$ and this must hold for any $M > 0$. This implies $\ssf(x^* + v) = c   (h(x) - h(x^*)) = 0$, which contradicts the assumption that $x^*$ is the unique global optimum. Hence, $C_u < \infty$.

The above argument proves that \ref{asm:radius} holds with the above constants for $a = 0$ and $b = \infty$. \Cref{prop:asm2} proves \ref{asm:levelset}.

\subsection{Proof of \Cref{pro:pupper}}\label{apdx:pro:pupper}
For a given $m \in \mathcal{X}_a^b$, there is a closed ball $\bar{\Ball}_u$ such that $S_0(m) \subseteq \bar{\Ball}_u$,
see \Cref{fig:levelset-balls}. We have
\begin{align*}
p^\mathrm{upper}_{(a, b]}(\bar\sigma)
&= \sup_{m \in \mathcal{X}_a^b} \sup_{\Sigma \in \mathcal{S}_\kappa} \int_{S_0(m)} \varphi\big(x; m, \left( \ssf(m) \bar\sigma \right)^2 \Sigma \big) dx \\
&\leq \sup_{m \in \mathcal{X}_a^b} \sup_{\Sigma \in \mathcal{S}_\kappa} \underbrace{ \int_{\bar{\Ball}_u} \varphi\big(x; m, \left( \ssf(m) \bar\sigma \right)^2 \Sigma \big) dx }_{(*1)}\enspace.
\end{align*}
The integral is maximized if the ball is centered at $m$. By a variable change ($x \leftarrow x - m$),
\begin{align*}
(*1) &\leq \int_{\norm{x} \leq C_u \ssf(m)} \varphi\big(x; 0, \left( \ssf(m) \bar\sigma \right)^2 \Sigma \big) dx
= \int_{\norm{x} \leq C_u / \bar\sigma} \varphi(x; 0, \Sigma) dx \\
&\leq \kappa^{d/2} \Phi\left(\bar{\Ball}\left(0, \frac{C_u}{\bar\sigma \kappa^{1/2}}\right); 0, \id \right)\enspace.
\end{align*}
Here we used $\Phi\big(\bar{\Ball}(0, r)\big); 0, \Sigma) \leq \kappa^{d/2} \Phi\left(\bar{\Ball}\left(0, \kappa^{-1/2} r \right); 0, \id \right)$ for any $r > 0$, which is proven in \Cref{lemma:phi} \new{below}. The right-most side (RMS) of the above inequality is independent of $m$. It proves \cref{eq:pupper}.

Similarly, there are balls $\Ball_{\ell}$ and $\bar{\Ball}_u$ such that $\Ball_{\ell} \subseteq S_0(m) \subseteq \bar{\Ball}_u$. We have
\begin{align*}
p^\mathrm{lower}_{(a, b]}(\bar\sigma)
&= \inf_{m \in \mathcal{X}_a^b} \inf_{\Sigma \in \mathcal{S}_\kappa} \int_{S_0(m)} \varphi\big(x; m,
\left( \ssf(m) \bar\sigma \right)^2 \Sigma \big) dx \\
&\geq \inf_{m \in \mathcal{X}_a^b} \inf_{\Sigma \in \mathcal{S}_\kappa} \underbrace{ \int_{\Ball_{\ell}} \varphi\big(x; m, \left( \ssf(m) \bar\sigma \right)^2 \Sigma \big) dx }_{(*2)}\enspace.
\end{align*}
The integral is minimized if the ball is at the opposite side of $m$ on the ball $\bar{\Ball}_u$,
see \Cref{fig:levelset-balls}.
By a variable change (moving $m$ to the origin) and letting $e_m = m / \norm{m}$,
\begin{align*}
(*2)
&
\geq \int_{\norm{x - ((2 C_u - C_{\ell}) \ssf(m)) e_m} \leq C_{\ell} \ssf(m)} \varphi\big(x; 0, \left( \ssf(m) \bar\sigma \right)^2 \Sigma\big) dx \\
&= \int_{\norm{x - ((2 C_u - C_{\ell}) / \bar\sigma)e_m} \leq C_{\ell} / \bar\sigma} \varphi(x; 0, \Sigma) dx \\
&\geq \kappa^{-d/2} \Phi\left(\bar{\Ball}\left(\left(\frac{(2 C_u - C_{\ell}) \kappa^{1/2}}{\bar\sigma}\right)e_m, \frac{C_{\ell} \kappa^{1/2}}{\bar\sigma} \right); 0, \id \right) \enspace.
\end{align*}
Here we used $\Phi\big( \bar\Ball(c, r); 0, \Sigma \big) \geq \kappa^{-d/2} \Phi\big(\bar\Ball(\kappa^{1/2} c, \kappa^{1/2} r); 0, \id \big)$ for any $c \in \R^d$ and $r > 0$ (\Cref{lemma:phi}). The RMS of the above inequality is independent of $m$ as its value is constant over all unit vectors $e_m$. Replacing $e_m$ with $e_1$, we have \cref{eq:plower}.

\new{
\begin{lemma}\label{lemma:phi}
  For all $\Sigma \in \mathcal{S}_\kappa$, 
  $\kappa^{-d/2} \varphi\left(x; 0, \kappa^{-1}\id \right) \leq \varphi\big(x; 0, \Sigma\big) \leq \kappa^{d/2} \varphi\left(x; 0, \kappa\id \right)$ and
  $\kappa^{-d/2} \Phi\left(\Ball(\sqrt{\kappa} c, \sqrt{\kappa} r); 0, \id \right) \leq \Phi\big(\Ball(c, r); 0, \Sigma\big) \leq \kappa^{d/2} \Phi\left(\Ball(c / \sqrt{\kappa}, r / \sqrt{\kappa}); 0, \id \right)$.
\end{lemma}
\begin{proof}
  For $\Sigma \in \mathcal{S}_\kappa$, we have $\det(\Sigma) = 1$ and $\Cond(\Sigma) = \lambda_{\max}(\Sigma) / \lambda_{\min}(\Sigma) \leq \kappa$.
  Since $\det(\Sigma) = 1$ and $\det(\Sigma) = \prod_{i=1}^{d}\lambda_i(\Sigma)$, we have $\lambda_{\max}(\Sigma) \geq 1 \geq  \lambda_{\min}(\Sigma)$.
  Therefore, we have $\lambda_{\min}(\Sigma) \geq \lambda_{\max} / \kappa \geq \kappa^{-1}$
  and $\lambda_{\max}(\Sigma) \leq \kappa \lambda_{\min}(\lambda) \leq \kappa$.
  Then we obtain $\kappa^{-1} x^\T\id x \leq  x^\T\Sigma^{-1} x \leq \kappa x^\T\id x$.
  With this inequality we have
  \begin{multline*}
    \varphi\big(x; 0, \Sigma\big)
    = (2\pi)^{-d/2} \exp(- x^\T \Sigma^{-1} x / 2)
    \leq (2\pi)^{-d/2} \exp(- x^\T \id x / (2\kappa)) \\
    = \kappa^{d/2} (2\pi\kappa)^{-d/2} \exp(- x^\T \id x / (2\kappa))
    = \kappa^{d/2}     \varphi\big(x; 0, \kappa \id\big) \enspace.
  \end{multline*}
  Analogously, we obtain $\varphi\big(x; 0, \Sigma\big) \geq \kappa^{-d/2} \varphi\big(x; 0, \kappa^{-1} \id\big)$.
  Taking the integral over $\Ball(c, r)$, we obtain the second statement.
\end{proof}%
}%

\subsection{Proof of \Cref{pro:exist_l_u}}\label{apdx:pro:exist_l_u}
  The upper bound of $p^\mathrm{upper}_{(a, b]}$ given in \eqref{eq:pupper} is strictly decreasing in $\bar\sigma$ and converges to zero when $\bar\sigma$ goes to infinity. This guarantees the
 existence of $\bar\sigma_u$ as a finite value. The existence of $\bar\sigma_{\ell} > 0$ is obvious under \ref{asm:levelset}.
	\ref{asm:radius} guarantees that there exists an open ball $B_{\ell}$ with radius $C_{\ell} (1-r)\ssf(m)$ such that $\Ball_{\ell} \subseteq \{x \in \R^d \mid \ssf(x) < (1-r)\ssf(m)\}$.
	Then, analogously to the proof of \Cref{pro:pupper}, the success probability with rate $r$ is lower bounded by
	\begin{equation}\textstyle
	 p^\mathrm{succ}_r(\bar\sigma; m, \Sigma) \geq \kappa^{-d/2} \Phi\left(\Ball\left(\left( \frac{(2 C_u - (1 - r) C_{\ell}) \kappa^{1/2}}{\bar\sigma} \right)e_1, \frac{(1-r) C_{\ell} \kappa^{1/2}}{\bar\sigma} \right); 0, \id\right) .
	 \end{equation}
	 The probability is independent of $m$, positive, and continuous in $\bar\sigma \in [\ell, u]$. Therefore the minimum is attained. This completes the proof.

\subsection{Proof of \Cref{pro:drift-UB-trunc}}\label{apdx:pro:drift-UB-trunc}

First, we remark that $m_t \in \mathcal{X}_{a, b}$ is equivalent to the condition $a < \ssf(m_t) \leq b$. If $\ssf(m_t) \leq a$ or $\ssf(m_t) > b$, \new{both sides of \cref{eq:Bdef} are zero}, hence the inequality is trivial. In the following we assume that $m_t \in \mathcal{X}_a^b$.

For the sake of simplicity we introduce
$\log^+(x) = \log(x)   \indlr{ x \geq 1 }$.
We rewrite the potential function as
\begin{align}\label{eq:a}
	V(\theta_t) = & \log \left( \ssf(m_t) \right)
	+ v   \log^+ \left( \frac{\aup  \elll  \ssf(m_t) }{\sigma_t}  \right)
	+ v   \log^+ \left( \frac{\sigma_t}{\adown  \uu  \ssf(m_t)} \right) \enspace.
\end{align}
The potential function at time $t+1$ can be written as
\begin{multline*}
V(\theta_{t+1}) = \log  \ssf(m_{t+1})
+ \underbrace{v   \log^{+}  \frac{\elll   \ssf(m_{t+1})}{\sigma_{t}}  \ind{ \sigma_{t+1} > \sigma_{t} }}_{P_2}
+ \underbrace{ v   \log^{+}  \frac{\aup   \elll   \ssf(m_{t})}{\adown   \sigma_{t}}  \ind{ \sigma_{t+1} < \sigma_{t} }}_{P_3} \\
		+ \underbrace{v   \log^{+}  \frac{\aup   \sigma_{t}}{\adown   \uu   \ssf(m_{t+1})}  \ind{ \sigma_{t+1} > \sigma_{t} }}_{P_4}
		+ \underbrace{v   \log^{+}  \frac{ \sigma_{t}}{\uu   \ssf(m_{t})}  \ind{ \sigma_{t+1} < \sigma_{t} }}_{P_5} \enspace.
\end{multline*}

We want to estimate the conditional expectation
\begin{equation}\label{eq:cond-expe}
\E\left[ \max\{ V(\theta_{t+1}) - V(\theta_{t}) \, ,\, -A \} \mid \theta_{t}\right].
\end{equation}

We partition the possible values of $\theta_t$ into three sets:
first the set of $\theta_t$ such that $\sigma_t < \elll   \ssf(m_t)$ ($\sigma_t$ is small),
second the set of $\theta_t$ such that $\sigma_t > \uu   \ssf(m_t)$ ($\sigma_t$ is large),
and last the set of $\theta_t$ such that $\elll   \ssf(m_t) \leq \sigma_t \leq \uu   \ssf(m_t)$ (reasonable $\sigma_t$).
In the following, we bound \cref{eq:cond-expe} for each of the
three cases and in the end our bound $B$ will equal the minimum
of the three bounds obtained for each case.

\textit{Reasonable $\sigma_t$ case: $\frac{\ssf(m_t)}{\sigma_t} \in \left[ \frac{1}{\uu}, \frac{1}{\elll} \right]$}. In case of success, where $\ind{\sigma_{t+1} > \sigma_t}=1$, we have $\ssf(m_{t+1}) / \sigma_{t+1} \leq \ssf(m_{t}) / (\aup \sigma_{t}) \leq 1 / (\aup \elll)$, implying that $P_2$ is always $0$. Similarly, in case of failure, $\ssf(m_{t+1}) / \sigma_{t+1} = \ssf(m_{t}) / (\adown \sigma_t) \geq 1 / (\adown \uu)$ and we find that $P_5$ is always zero. We rearrange $P_3$ and $P_4$ into
\begin{align*}
	P_3 &= v   \log^+ \left( \frac{\aup   \elll   \ssf(m_{t})}{\adown   \sigma_{t}} \right) \ind{ \sigma_{t+1} < \sigma_t } \enspace, \\
	P_4 &= v  \left[ \log \left( \frac{\aup   \sigma_{t}}{\adown   \uu   \ssf(m_{t})} \right) - \log \left( \frac{ \ssf(m_{t+1}) }{ \ssf(m_{t}) } \right)  \right] \ind{\frac{\adown \uu \ssf(m_{t+1})}{\aup \sigma_{t}} < 1 } \ind{ \sigma_{t+1} > \sigma_t } \enspace.
\end{align*}
Then, the one-step change $\Delta_{t} = V(\theta_{t+1}) - V(\theta_{t})$ is upper bounded by
\begin{multline}
	\Delta_{t} \leq  \left(1 - v   \ind{\frac{\adown \uu \ssf(m_{t})}{\aup \sigma_{t}} < 1 }   \ind{ \sigma_{t+1} > \sigma_t } \right) \log \left( \frac{ \ssf(m_{t+1}) }{ \ssf(m_{t}) } \right)
	\\
	 + v   \log^+ \left( \frac{\aup   \elll   \ssf(m_{t})}{\adown   \sigma_{t}} \right)    \ind{ \sigma_{t+1} < \sigma_t } + v   \log^+ \left( \frac{\aup   \sigma_{t}}{\adown   \uu   \ssf(m_{t})} \right)   \ind{ \sigma_{t+1} > \sigma_t }  \\
	\leq  (1 - v )  \log  \frac{ \ssf(m_{t+1}) }{ \ssf(m_{t}) }
 + v  \log^+  \frac{\aup  \elll  \ssf(m_{t})}{\adown  \sigma_{t}}    \ind{ \sigma_{t+1} < \sigma_t }
	 + v  \log^+  \frac{\aup  \sigma_{t}}{\adown  \uu  \ssf(m_{t})}   \ind{ \sigma_{t+1} > \sigma_t }
	\enspace.
\end{multline}
The truncated one-step change $\max\{ \Delta_{t} \, , \, - A\}$ is upper bounded by
\begin{multline}
	\max\{ \Delta_{t} \, , \, - A\}
	\leq (1 - v ) \max\left\{  \log \left( \frac{ \ssf(m_{t+1}) }{ \ssf(m_{t}) } \right)  \, , \, - \frac{A}{1 - v} \right\}\\
	+ v   \log^+ \left( \frac{\aup   \elll   \ssf(m_{t})}{\adown   \sigma_{t}} \right)    \ind{ \sigma_{t+1} < \sigma_t } + v  \log^+ \left( \frac{\aup  \sigma_{t}}{\adown  \uu  \ssf(m_{t})} \right)  \ind{ \sigma_{t+1} > \sigma_t }
	\enspace.
\end{multline}
To consider the expectation of the above upper bound, we need to compute
the expectation of the maximum of
$\log \left( \frac{ \ssf(m_{t+1}) }{ \ssf(m_{t}) } \right)$ and
$-\frac{A}{1 - v}$. Let $a \leq 0$ and $b \in \R$ then $\max(a, b) = a   \ind{a > b} + b   \ind{a \leq b} \leq b   \ind{a \leq b}$. Applying this and taking the conditional expectation, a trivial upper
bound for the conditional expectation of
$%
\max\left\{ \log \left( \frac{ \ssf(m_{t+1}) }{ \ssf(m_{t}) } \right)  \, , \, - \frac{A}{1 - v} \right\}
$
is $-\frac{A}{1 - v}$ times the probability of
$\log \left( \frac{ \ssf(m_{t+1}) }{ \ssf(m_{t}) } \right)$
being no greater than $-\frac{A}{1 - v}$.
The latter condition is equivalent to
$\ssf(m_{t+1}) \leq (1-r)   \ssf(m_{t})$ corresponding to successes with
rate $r = 1 - \exp\left( - \frac{A}{1 - v}\right)$ or better.
That is,
\begin{equation}\textstyle
	(1 - v )  \E \left[ \max\left\{  \log \left( \frac{ \ssf(m_{t+1}) }{ \ssf(m_{t}) } \right)  \, , \, - \frac{A}{1 - v} \right\} \right] \\
	\leq - A  p^\mathrm{succ}_{r}\left(\frac{ \sigma_t}{\ssf(m_t)}; m_t, \Sigma_t\right)
	\enspace.
\end{equation}
Note also that the expected value of $\indlr{ \sigma_{t+1} > \sigma_t }$
is the success probability, namely, $p^\mathrm{succ}_{0}\left(\frac{\sigma_t}{\ssf(m_t)}; m_t, \Sigma_t\right)$.
We obtain an upper bound for the conditional expectation of
$\max\{ \Delta_{t} \, , \, - A\}$ in the case of reasonable $\sigma_t$
as
\begin{multline}
	\E\left[\max\{ \Delta_{t} \, , \, - A\} | \theta_t \right] \leq - A  p^\mathrm{succ}_{r}\left(\frac{\sigma_t}{\ssf(m_t)}; m_t, \Sigma_t\right)   \\
	+ \bigg( \log \left(\frac{\aup}{\adown}\right) + \underbrace{ \log\left( \frac{ \ell   \ssf(m_{t})}{ \sigma_t} \right)}_{\leq 0} \bigg)  v  \left( 1 - p^\mathrm{succ}_{0}\left(\frac{ \sigma_t}{\ssf(m_t)}; m_t, \Sigma_t\right) \right) \\
	+ \bigg( \log \left(\frac{\aup}{\adown}\right) + \underbrace{\log \left( \frac{ \sigma_t}{\uu  \ssf(m_{t})} \right)}_{\leq 0} \bigg)   v  p^\mathrm{succ}_{0}\left(\frac{ \sigma_t}{\ssf(m_t)}; m_t, \Sigma_t \right)
	\leq - A  p^*_r + v  \log \left(\frac{\aup}{\adown}\right)\enspace.
	\label{eq:case3bound}
\end{multline}

\textit{Small $\sigma_t$ case: $\frac{\ssf(m_t)}{ \sigma_t} > \frac{1}{\elll}$}. If $\ell \ssf(m_t) > \sigma_t$, the 2nd summand in \cref{eq:a} is positive. Moreover, if $\sigma_{t + 1} < \sigma_{t}$, we have
$\elll \ssf(m_{t+1}) = \elll \ssf(m_{t}) > \sigma_t > \sigma_{t+1}$ and hence the 2nd summand in \cref{eq:a} is positive for $V(\theta_{t+1})$ as well. If $\sigma_{t + 1} > \sigma_{t}$, any regime can happen. Then, $V(\theta_{t+1}) - V(\theta_{t}) =$
\begin{align*}
	=& \log  \frac{\ssf(m_{t+1})}{ \ssf(m_{t}) }  - v  \log  \frac{\aup  \elll  \ssf(m_t)}{\sigma_t} \notag + v  \log  \frac{\elll  \ssf(m_{t+1})}{\sigma_{t}}  \indlr{ \frac{\elll \ssf(m_{t+1})}{\sigma_{t}} > 1 }\indlr{ \sigma_{t+1} > \sigma_{t} } \notag \\
         &+ v  \log  \frac{\aup  \elll  \ssf(m_{t})}{\adown  \sigma_{t}}  \indlr{ \frac{\aup \elll \ssf(m_{t})}{\adown \sigma_{t}} > 1 }\indlr{ \sigma_{t+1} < \sigma_{t} } \notag \\
  &+ v  \log  \frac{\aup  \sigma_{t}}{\adown  \uu  \ssf(m_{t+1})}  \indlr{\frac{\adown \uu \ssf(m_{t+1})}{\aup \sigma_{t}} < 1 } \indlr{ \sigma_{t+1} > \sigma_{t} } \notag \\
	=& \log \left( \frac{\ssf(m_{t+1})}{ \ssf(m_{t}) } \right) \left[1 + v   \left(\indlr{ \frac{\elll \ssf(m_{t+1})}{\sigma_{t}} > 1 } - \indlr{\frac{\adown \uu \ssf(m_{t+1})}{\aup \sigma_{t}} < 1 }  \right)   \indlr{ \sigma_{t+1} > \sigma_{t} } \right] \notag \\
	&- v   \log \left( \frac{\adown   \uu   \ssf(m_{t})}{\aup   \sigma_{t}} \right) \indlr{\frac{\adown \uu \ssf(m_{t+1})}{\aup \sigma_{t}} < 1 } \indlr{ \sigma_{t+1} > \sigma_{t} } \notag \\
	&- v   \log \left( \frac{\elll   \ssf(m_{t})}{\sigma_t} \right) \left[ 1 - \indlr{ \frac{\elll \ssf(m_{t+1})}{\sigma_{t}} > 1 }\indlr{ \sigma_{t+1} > \sigma_{t} } - \indlr{ \frac{\aup \elll \ssf(m_{t})}{\adown \sigma_{t}} > 1 }\indlr{ \sigma_{t+1} < \sigma_{t} }  \right] \notag \\
	&- v   \left( \log (\aup) - \log \left(\frac{\aup}{\adown}\right) \indlr{ \frac{\aup \elll \ssf(m_{t})}{\adown \sigma_{t}} > 1 }\indlr{ \sigma_{t+1} < \sigma_{t} } \right) \enspace.
    \label{eq:case1}
\end{align*}
On the RMS of the above equality, the first term is guaranteed to be
non-positive since $v \in (0, 1)$.
The second and third terms are non-positive as well since
$\frac{\adown \uu \ssf(m_{t})}{\aup \sigma_{t}} > \frac{\adown\uu}{\aup \elll} \new{\geq} 1$
and $\frac{\elll \ssf(m_{t})}{\sigma_{t}} > 1$. Replacing the indicator $\indlr{ \frac{\aup \elll \ssf(m_{t})}{\adown \sigma_{t}} > 1 }$ with $1$ in the last term provides an upper bound. Altogether, we obtain
\begin{equation*}
	\Delta_t = V(\theta_{t+1}) - V(\theta_{t})
	\leq - v    \left( \log ( \aup ) - \log ( \aup / \adown ) \indlr{ \sigma_{t+1} < \sigma_{t} } \right)  \enspace.
\end{equation*}
Note that the RHS is larger than $- A$ since it is lower bounded by $- v   \log(\aup)$ and $v \leq A / \log(\aup)$. Then, the conditional
expectation of $\max\{ \Delta_{t} \, , \, - A\}$ is
\begin{multline}
  \E\left[\max\{ \Delta_{t} \, , \, - A\} | \mathcal{F}_t \right]
   \leq - v    \left( \log\left( \frac{\aup}{\adown} \right) p^\mathrm{succ}_{0}\left( \frac{\sigma_t}{\ssf(m_t)}; m_t, \Sigma_t \right) + \log( \adown )  \right) \\
   \leq - v    \left( \log\left( \frac{\aup}{\adown} \right) p_{\ell} + \log( \adown )  \right)
   \new{= - v    \log\left( \frac{\aup}{ \adown} \right) \left( p_{\ell} - p_\mathrm{target}  \right) }
   = - v    \frac{p_{\ell} - p_u}{2}   \log\left( \frac{\aup}{\adown} \right) \enspace.
 \label{eq:case1bound}
\end{multline}
Here we used $\mathbb{E}[1\{\sigma_{t+1}<\sigma_t\}\mid\mathcal{F}_t] = 1-p_0^{\mathrm{succ}}\left(\frac{\sigma_t}{f_\mu(m_t)};m_t, \Sigma_t\right)$ for the first inequality, 
$p^\mathrm{succ}_{0}\left( \frac{\sigma_t}{\ssf(m_t)}; m_t, \Sigma_t \right) > p_\ell$ for the second inequality,
$p_\mathrm{target} = \log\left(\frac{1}{\adown}\right) \Big/ \log\left(\frac{\aup}{\adown}\right)$ and $p_\mathrm{target} = (p_\uu+p_\elll)/2$ for the last two equalities.

\textit{Large $\sigma_t$ case: $\frac{\ssf(m_t)}{\sigma_t} < \frac{1}{\uu}$}.
Since $\frac{\ssf(m_{t+1})}{\sigma_{t+1}} \leq \frac{\ssf(m_t)}{\adown \sigma_t} < \frac{1}{\adown \uu}$,
the 3rd summand in \cref{eq:a} is positive in both $V(\theta_t)$ and $V(\theta_{t+1})$.
For the 2nd summand in \cref{eq:a}, recall that $\aup \elll \ssf(m_{t})/ \sigma_t < \aup \elll / \uu \leq \adown < 1$
since we have assumed that $\uu / \elll \geq \aup / \adown$. Hence, for
$V(\theta_t)$ the 2nd summand in \cref{eq:a} is zero. Also,
$\aup \elll \| m_{t+1} \| / \sigma_{t+1} \leq \aup \elll / (\adown \uu ) = (\aup / \adown) \elll / \uu \geq 1$
and thus for $V(\theta_{t+1})$ the 2nd summand in \cref{eq:a} also equals $0$.
We obtain
\begin{equation*}
	V(\theta_{t+1}) - V(\theta_t)
	=  (1 - v) \big( \log\left( \ssf(m_{t+1}) \right) - \log\left( \ssf(m_{t}) \right) \big)
	   + v   \log\left( \sigma_{t+1} / \sigma_t \right).
\end{equation*}
The first term on the RHS is guaranteed to be non-positive since $v < 1$,
yielding $\Delta_{t} \leq v   \log(\sigma_{t+1}/\sigma_t)$.
On the other hand,
\begin{align*}
	v   \log(\sigma_{t+1}/\sigma_t) 
	& =v   \left( \log (\aup)  \ind{ \sigma_{t+1} > \sigma_{t} } + \log(\adown) \ind{ \sigma_{t+1} < \sigma_{t} } \right) \\
	& =v   \left( \log(\aup / \adown) \ind{ \sigma_{t+1} > \sigma_{t} } - \log(1/\adown) \right) \\
	& \geq - v \log(1/\adown) \geq -A \enspace,
\end{align*}
where the last inequality comes from the prerequisite $v \leq A / \log(1 / \adown)$. 
Hence,
\begin{equation*}
\max\{ \Delta_{t} \, ,\allowbreak \, - A\} \leq \max \{ v   \log(\sigma_{t+1}/\sigma_t) , - A\} = v \log(\sigma_{t+1} / \sigma_t) \enspace.
\end{equation*}
Then, the conditional expectation of $\max\{ \Delta_{t} \, ,\allowbreak \, - A\}$ is
\begin{multline}
	\E\left[\max\{ \Delta_{t} \, , \, - A\} | \theta_t \right]
	 \leq v \left(\log(\adown) + \log\left(\frac{\aup}{\adown}\right) p^\mathrm{succ}_{0}\left( \frac{\sigma_t}{\ssf(m_t)}; m_t, \Sigma_t \right) \right) \\
	 \leq v \left(\log(\adown) + \log\left(\frac{\aup}{\adown}\right) p_{u} \right)
	\new{= v \log\left(\frac{\aup}{\adown}\right)\left(-p_\mathrm{target} +  p_{u} \right) }
    = - v    \frac{p_{\ell} - p_u}{2}   \log\left( \frac{\aup}{\adown} \right) \enspace.
	\label{eq:case2bound}
\end{multline}
Here we used $p^\mathrm{succ}_{0}\left( \frac{\sigma_t}{\ssf(m_t)}; m_t, \Sigma_t \right) \leq p_{u}$.

\textit{Conclusion}.
Inequalities \cref{eq:case1bound,eq:case2bound,eq:case3bound} together cover all possible cases and we hence obtain \cref{eq:B}.

Finally, we prove the positivity of $B$ for an arbitrary $A > 0$. \Cref{pro:exist_l_u} guarantees the positivity of $p_r^*$ for any choice of $A$ since $r = 1 - \exp(- A / (1 - v)) \in \left(0, 1\right)$ for any $A > 0$ and $v < 1$. Therefore, $A   p_r^* > 0$ for any $A$ and $v \leq \min(1,\ A/\log(1/\adown),\ A / \log(\aup))$. Moreover, for a sufficiently small $v$, $p^*_r$ is strictly positive for any $A > 0$. Therefore, one can take a sufficiently small $v$ that satisfies $A   p_r^* > v \log(\aup / \adown)$. The first term in the minimum in \cref{eq:B} is positive. The second term therein is clearly positive for $v > 0$. This completes the proof.

\subsection{Proof of \Cref{proposition:scaling}}\label{apdx:proposition:scaling}

Consider $d \geq 2$. We set $A = 1/d$. We bound $B$ from below by taking a specific value for $v \in (0,\ \min(1, A / \log(1/\adown),\ A / \log(\aup))$ instead of considering $\sup$ for $v$.
Our candidate is $v = \frac{A   p'}{\log(\aup / \adown)}   \frac{2}{(2 + p_{\ell} - p_u)}$, where
$p' = \inf_{\barsigma \in [\elll, \uu]} p_{r'}(\barsigma)$
and
$r' = 1 - \exp\big(-A   \big(1 - \frac{1}{d   \log(\aup/\adown)}\big)^{-1}\big)$.
It holds $v < \frac{1}{d   \log(\aup/\adown)}$ and hence
$r' > r$, from which we obtain $p' < p^*$.

We bound the terms in \cref{eq:B} as:
   $A p^* -  v \log(\aup / \adown) = \frac{p'}{d}\left(\frac{p^*}{p'} - \frac{2}{2 + p_{\ell} - p_u}\right) \geq \frac{p'}{d}\left(\frac{p_{\ell} - p_u}{2 + p_{\ell} - p_u}\right)$
   and
   $v \frac{p_{\ell} - p_u}{2} \log\left(\frac{\aup}{\adown}\right) = \frac{p'}{d} \frac{p_{\ell} - p_u}{2 + p_{\ell} - p_u}$.
Therefore, we have $B \geq \frac{p'}{d} \frac{p_{\ell} - p_u}{2 + p_{\ell} - p_u}$.
Note that one can take $p_\ell - p_u \in \Theta(1)$ since the only condition is $p_\mathrm{target} = (p_\ell + p_u) / 2 \in \Theta(1)$. To obtain $B \in \Omega(1/d)$, it is sufficient to show $p' \in \Theta(1)$ for $d \to \infty$.

  Fix $p_\ell$ and $p_u$ independently of $d$.
  In the light of Lemma~3.1 in \cite{akimoto2018drift}, we have that $p_0: \R_{>} \to (0, 1/2)$ is continuous and strictly decreasing from $1/2$ to $0$ for all $d \in \mathbb{N}$.
  Therefore, for each $d \in \mathbb{N}$ there exists an inverse map $p_0^{-1}: (0, 1/2) \to \R_{>}$.
  Define $\hat\sigma_\ell^{d} = d   V_d   p_0^{-1}(p_\ell)$ and $\hat\sigma_u^{d} = d   V_d   p_0^{-1}(p_u)$ for each $d \in \mathbb{N}$.
  It follows from Lemma~3.2 in \cite{akimoto2018drift} that \new{$p_0^{\mathrm{lim}}:\bar\sigma\mapsto\lim_{d\to\infty}p_0(\bar\sigma)$} is also strictly decreasing, hence invertible.
  \new{The existence of $\lim_{d\to\infty}p_0(\cdot)$ is also proved in \cite{akimoto2018drift}}.
  We let $\hat\sigma_\ell^{\infty} = (p_0^{\mathrm{lim}})^{-1}(p_\ell)$ and $\hat\sigma_u^{\infty} = (p_0^{\mathrm{lim}})^{-1}(p_u)$.
  Because of the pointwise convergence of $p_0(\barsigma = \hat\sigma / (d V_d))$ to $p_0^\mathrm{lim}(\hat\sigma)$,
  we have $\hat\sigma_\ell^{d} \to \hat\sigma_\ell^{\infty}$ and $\hat\sigma_u^{d} \to \hat\sigma_u^{\infty}$ for $d \to \infty$.
  Hence, for any $\hat{u} > \hat\sigma_u^{\infty}$ and $\hat{\ell} < \hat\sigma_\ell^{\infty}$ with $u / \ell \geq \aup / \adown$,
  there exists $D \in \mathbb{N}$ such that for all $d \geq D$ we have $\hat{u} > \hat\sigma_u^{d}$ and $\hat{\ell} < \hat\sigma_\ell^{d}$.
  Now we fix $\hat{u}$ and $\hat{\ell}$ in this way. This amounts to selecting $u = d   V_d   \hat{u}$ and $\ell = d   V_d   \hat{\ell}$.


We have $\lim_{d \to \infty} d   r' = 1$ since $\lim_{d \to \infty} d   \log(\aup / \adown) = \infty$ and hence according to
Lemma~3.2 in \cite{akimoto2018drift} we have
\begin{align*}
\liminf_{d \to \infty} p'
 		&= \liminf_{d \to \infty} \min_{\barsigma \in [\elll, \uu]} \left\{ p_{r'} (\barsigma) \right\}
 		= \liminf_{d \to \infty} \min_{\hat\sigma \in [\hat{\elll}, \hat{\uu}]} p_{r'} \left(\frac{\hat\sigma}{d   V_d}\right)
 		\\
 		&\overset{(\star)}{=} \min_{\hat\sigma \in [\hat\elll, \hat\uu]} \lim_{d \to \infty} \left( p_{r'} \left(\frac{\hat\sigma}{d   V_d}\right)\right)
 		= \min_{\hat\sigma \in [\hat\elll, \hat\uu]} \Psi\left(-\frac{1}{\hat\sigma} - \frac{{\hat\sigma}}{2}\right)\enspace,
\end{align*}
where the equality $(\star)$ follows from the pointwise convergence of $p_{r'}$ to $\lim_{d \to \infty} p_{r'}$ and the continuity of $p_{r'}$ and $\lim_{d \to \infty} p_{r'}$.%
\footnote{%
Let $\{f_n : n \geq 1\}$ be a sequence of continuous functions on $\R$ and $f$ be a continuous function such that $f$ is the pointwise limit $\lim_{n} f_n(x) = f(x)$ of the sequence. Since they are continuous, there exist the minimizers of $f_n$ and $f$ in a compact set $[\ell, u]$. Let $x_n = \argmin f_n(x)$ and $x^* = \argmin f(x)$, where $\argmin$ is taken over $x \in [\ell, u]$ and we pick one if there exist more than one minimizers. It is easy to see that $f_n(x_n) \leq f_n(x^*)$, hence $\liminf_n f_n(x_n) \leq \liminf_{n} f_n (x^*) = f(x^*)$. Let $\{n_i : i \geq 1\}$ be the sub-sequence of the indices such that $\liminf_{n} f_n(x_n) = \lim_{i} f_{n_i}(x_{n_i})$. Since $\{x_{n_i} : i \geq 1\}$ is a bounded sequence, Bolzano-Weirstra{\ss} theorem provides a convergent sub-sequence $\{x_{n_{i_k}} : k \geq 1\}$ and we denote its limit as $x_*$. Of course we have $\liminf_{n} f_n(x_n) = \lim_{k} f_{n_{i_k}}(x_{n_{i_k}})$. Due to the continuity of $\{f_n : n \geq 1\}$ and the pointwise convergence to $f$, we have $\lim_{k} f_{n_{i_k}}(x_{n_{i_k}}) = \lim_{k} f_{n_{i_k}}(x_*) = f(x_*)$. Therefore, $\liminf_{n} f_n(x_n) = f(x_*) \leq f(x^*)$. Since $x^*$ is the minimizer of $f$ in $[\ell, u]$ and $x_* \in [\ell, u]$, it must hold $f(x_*) \geq f(x^*)$. Hence, $\liminf_{n} f_n(x_n) = f(x^*)$.}
This completes the proof.

\end{document}